\documentclass{amsart}
\usepackage{amsmath,amsthm,amssymb}
\usepackage{graphicx}

\newtheorem{theorem}{Theorem}[section]
\newtheorem{proposition}[theorem]{Proposition}
\newtheorem{lemma}[theorem]{Lemma}
\newtheorem{corollary}[theorem]{Corollary}

\theoremstyle{definition}
\newtheorem{definition}[theorem]{Definition}
\newtheorem{remark}[theorem]{Remark}

\newtheorem{example}[theorem]{Example}
\newtheorem{question}[theorem]{Question}

\title[Nuclei and exotic 4-manifolds]{Nuclei and exotic 4-manifolds}
 
\author[Kouichi Yasui]{Kouichi Yasui}
\thanks{The author was partially supported by KAKENHI 23840027.}
\date{February 17, 2012}
\subjclass[2010]{Primary~57R55, Secondary~57R65, 57R17}
\keywords{4-manifold; smooth structure; cork; Stein manifold; 3-manifold}
\address{Department~of~Mathematics, Graduate School~of~Science, Hiroshima~University, 1-3-1 Kagamiyama, Higashi-Hiroshima, 739-8526, Japan}\email{kyasui@hiroshima-u.ac.jp}

\begin{document}

\begin{abstract}We introduce a new generalization of Gompf nuclei and give applications. We construct infinitely many exotic smooth structures for a large class of compact 4-manifolds with boundary, regarding topological invariants. We prove that a large class of closed 3-manifolds (including disjoint unions of Stein fillable $3$-manifolds) bound compact connected oriented 4-manifolds which admit infinitely many smooth structures. To detect smooth structures, we introduce a relative genus function. 
As a side result, we show that log transform and knot surgery of 4-manifolds do not admit any Stein structure, under a mild condition. Applying these results together with corks, we construct arbitrary many compact Stein 4-manifolds and infinitely many non-Stein 4-manifolds which are all homeomorphic but mutually non-diffeomorphic. 
\end{abstract}

\maketitle

\section{Introduction}\label{sec:intro}

It is currently known that many closed (or non-compact) 4-manifolds admit infinitely many smooth structures (cf.\ \cite{GS}, \cite{S}, \cite{FS4}). However, smooth structures on compact 4-manifolds with boundary are not much investigated, partially due to the difficulty of the invariants (e.g.\ relative Seiberg-Witten invariants). The main purpose of this paper is to study smooth and Stein structures on such 4-manifolds. 

Gompf~\cite{G0} (and independently Ue~\cite{Ue}) introduced nuclei, which is useful to construct exotic smooth structures of ``closed'' 4-manifolds (e.g.\ \cite{BDPS}, \cite{P2}, \cite{APU}). They also proved that each Gompf nucleus admits infinitely many smooth structures. 

In this paper, we introduce a new generalization of Gompf nuclei. 
Unlike Fuller's generalization for higher genus fibrations in \cite{Fu}, we require our nucleus to contain a cusp neighborhood. 
See Definition~\ref{def:nucleus} for details. Applying nuclei, we construct exotic smooth structures for a large class of non-closed 4-manifolds, regarding topological invariants and boundary 3-manifolds. 

We first prove the theorem below using log transform and knot surgery, where the closed case follows from theorems of Boyer~\cite{B} and Fintushel-Stern~\cite{FS1}, \cite{FS2}. 

\begin{theorem}\label{sec:intro:thm:original} Let $X$ be an arbitrary connected oriented smooth $4$-manifold. Here $X$ is allowed to be closed, or to have $($possibly disconnected$)$ boundary, or to be non-compact. Suppose that $X$ smoothly contains a nucleus. Assume further that $X$ can be embedded into a connected oriented closed smooth $4$-manifold $Z$ with $b^+_2>1$ whose Seiberg-Witten invariant does not vanish. Then there exist infinitely many smooth 4-manifolds which are all homeomorphic to $X$ but mutually non-diffeomorphic. 
\end{theorem}
\begin{remark}This theorem also holds in the case where the closed 4-manifold $Z\supset X$ satisfies $\pi_1(Z)\cong 1$, $b_2^+(Z)=1$, $b_2^-(Z)\leq 9$ and $SW_Z\neq 0$ (see Subsection~\ref{subsection:remarks}). 
\end{remark}

To detect smooth structures, we introduce a relative genus function (Definition \ref{def:relative genus}), generalizing techniques of Akbulut and the author in~\cite{AY5}. 
While the reader may expect that the relative Seiberg-Witten (or Ozsv\'{a}th-Szab\'{o}) invariants also detect smooth structures in the compact case, it seems that known formulas (to the author) do not work in this general setting (cf.\ Theorem~3.3 of \cite{Mar}). 

In the rest of this section, except for Theorem~\ref{intro:non-stein}, we discuss applications of Theorem~\ref{sec:intro:thm:original}. Here recall that a ($4$-dimensional oriented)
$2$-handlebody means a compact, connected, oriented smooth $4$-manifold obtained from the $4$-ball by attaching $1$- and $2$-handles. Applying Theorem~\ref{sec:intro:thm:original} and corks, we can easily construct infinitely many exotic smooth structures for a large class of $4$-manifolds with boundary, regarding topological invariants. 

\begin{theorem}\label{sec:intro:thm:modified}
Let $X$ be a $2$-handlebody which contains a nucleus as a subhandlebody. 
Suppose that the nucleus admits a Stein structure. 
Then, there exists a compact connected oriented smooth $4$-manifold $X_0$ which satisfies the following.
\\
$(1)$ There exist infinitely many smooth 4-manifolds which are all homeomorphic to $X_0$ but mutually non-diffeomorphic.
\\
$(2)$ The fundamental group, the integral homology groups, the integral homology groups of the boundary, and the intersection form of $X_0$ are isomorphic to those of $X$. 
\\
$(3)$ $X_0$ $($resp.\ $X$$)$ can be embedded into $X$ $($resp.\ $X_0$$)$. 
\end{theorem}

See also Corollaries~\ref{cor:exotic:2-handlebody:boundary sum} and \ref{sec:infinite exotic:cor} for simple consequences. Note that Theorem~1.1 of \cite{AY5} gives finitely many exotic smooth structures for a larger class. 


We next consider smooth structures of $4$-manifolds from a view point of boundary $3$-manifolds. Since every closed connected oriented $3$-manifold bounds a simply connected compact oriented smooth $4$-manifold, the following question is natural. 
\begin{question}
Does every closed oriented (possibly disconnected) $3$-manifold bounds a (simply) connected compact oriented $4$-manifold which has infinitely many distinct smooth structures?
\end{question}There are a few affirmative examples. 
$\#_n(S^1\times S^2)$ $(n\geq 0)$ is easy to check using closed elliptic surfaces. Other known examples (to the author) are certain homology $3$-spheres (Gompf~\cite{G0} and Ue~\cite{Ue}), certain circle bundles over surfaces (Fintushel-Stern~\cite{FS1.6}, \cite{FS1.5} and Mark~\cite{Mar}) and certain Seifert 3-manifolds (Akhmedov-Etnyre-Mark-Smith~\cite{AEMS}). 
In this paper, we give an affirmative answer for a large class of $3$-manifolds applying Theorem~\ref{sec:intro:thm:original}.
 Note that most connected closed orientable $3$-manifolds are known to be Stein fillable (\cite{G}, \cite{GS}, \cite{OS1}), though there are infinitely many non Stein fillable $3$-manifolds (\cite{L}, \cite{LS1}, \cite{LS2}). 
\begin{theorem}\label{intro:Stein boundary}$(1)$ For every Stein fillable $3$-manifold $M$, there exist infinitely many smooth 4-manifolds $X_n$ $(n\geq 1)$ such that they are all homeomorphic but mutually non-diffeomorphic and that the boundary $\partial X_n$ of each $X_n$ is diffeomorphic to $M$. 
\smallskip\\
$(2)$ Let $M$ be a disjoint union of arbitrary finite number of Stein fillable $3$-manifolds. Then there exist infinitely many smooth 4-manifolds $X_n$ $(n\geq 1)$ such that they are all homeomorphic but mutually non-diffeomorphic and that the boundary $\partial X_n$ of each $X_n$ is diffeomorphic to $M$.
\end{theorem}

Moreover, ignoring the simple connectivity of $4$-manifolds, we also prove this for more general $3$-manifolds including non Stein fillable $3$-manifolds with both orientations. See Theorem~\ref{thm:boundary:symplectic}, Corollary~\ref{cor:symplectic:boundary} and Example~\ref{ex:symplectic:boundary}. 


We here discuss non-existence of Stein structures. It is interesting to find surgical operations which produce compact Stein $4$-manifolds. Since log transform (resp.\ knot surgery) produces closed complex surfaces (resp.\  symplectic $4$-manifolds) under certain conditions (cf.\ \cite{GS}), one might expect that these operations can produce compact Stein $4$-manifolds. However, we give a negative answer under a mild condition. 

\begin{theorem}\label{intro:non-stein}Let $X$ be a compact connected oriented smooth $4$-manifold with boundary. Suppose that $X$ contains a $c$-embedded torus $T$ and that $T$ represents a non-torsion $($i.e.\ infinite order$)$ class  of $H_2(X;\mathbb{Z})$. Then the following hold. 
\\
$(1)$ For each $p\geq 2$, the $p$-log transform $X_{(p)}$ of $X$ along $T$ does not admit any Stein structure for both orientations.
\\
$(2)$ For each knot $K$ in $S^3$ with the non-trivial Alexander polynomial, the knot surgery $X_K$ of $X$ along $T$ does not admit any Stein structure for both orientations.
\end{theorem}

We finally discuss smooth structures of Stein $4$-manifolds. It is known that diffeomorphism types of Stein $4$-manifolds bounded by certain $3$-manifolds are unique $($e.g.\ $\#_n\,S^1\times S^2$ $(n\geq 0)$. cf.\ \cite{OS1}). 
It is thus interesting to find exotic Stein $4$-manifold pairs. Akhmedov-Etnyre-Mark-Smith~\cite{AEMS} constructed the first example of infinitely many compact Stein $4$-manifolds which are all homeomorphic but mutually non-diffeomorphic. Akbulut and the author (Theorem~10.1 of \cite{AY5}) recently constructed arbitrary many compact Stein $4$-manifolds and arbitrary many non-Stein $4$-manifolds which are all homeomorphic but mutually non-diffeomorphic. In this paper, we construct arbitrary many compact Stein $4$-manifolds and infinitely many non-Stein $4$-manifolds which are all homeomorphic but mutually non-diffeomorphic.

\begin{theorem}\label{intro:stein and non-stein}
Let $X$ be a $2$-handlebody which contains a nucleus as a subhandlebody. 
Suppose that the nucleus admits a Stein structure. Then for each integer $n\geq 1$, there exist infinitely many compact connected oriented smooth $4$-manifolds $X_i$ $(i=0,1,2,\dots)$ which satisfy the following. 
\\
$(1)$ $X_i$ $(i\geq 0)$ are all homeomorphic but mutually non-diffeomorphic. 
\\
$(2)$ Each $X_i$ $(1\leq i\leq n)$ admits a Stein structure, and any $X_i$ $(i\geq n+1)$ admits no Stein structure. 
\\
$(3)$ The fundamental group, the integral homology groups, the integral homology groups of the boundary, and the intersection form of each $X_i$ $(i\geq 0)$ are isomorphic to those of $X$.
\\
$(4)$ Each $X_i$ $(0\leq i\leq n)$ can be embedded into $X$ . 
\\
$(5)$ $X$ can be embedded into $X_0$. 
\end{theorem}
See also Corollaries~\ref{cor:nonstein:2-handlebody:boundary sum} and \ref{cor:nonstein:group} for simple consequences. In a forthcoming paper~\cite{Y5}, we will apply Theorem~\ref{sec:intro:thm:original} to construct exotic $S^2$-knots and links in certain 4-manifolds.

This paper is organized as follows. In Section~\ref{sec:log}, we recall log transform, knot surgery and the adjunction inequality. In Section~\ref{sec:nucleus}, we study a generalization of Gompf nuclei. In Section~\ref{sec:minimal genus}, we introduce a relative genus function. In Section~\ref{sec:nuclei and exotic}, we prove Theorem~\ref{sec:intro:thm:original}. In Sections~\ref{sec:stein} and \ref{sec:cork}, we briefly review compact Stein $4$-manifolds, corks and $W$-modifications. In Section~\ref{sec:stein and exotic}, we prove Theorem~\ref{sec:intro:thm:modified}, applying $W$-modifications and Theorem~\ref{sec:intro:thm:original}. In Section~\ref{sec:boundary}, we prove Theorem~\ref{intro:Stein boundary}, using compact Stein $4$-manifolds and applying Theorem~\ref{sec:intro:thm:original}. In Section~\ref{sec:non-existence}, we prove Theorem~\ref{intro:non-stein}. In Section~\ref{sec:stein and non-stein}, we recall an algorithm of Akbulut and the author~\cite{AY5} and prove Theorem~\ref{intro:stein and non-stein}, applying the algorithm and Theorems~\ref{sec:intro:thm:original} and \ref{intro:non-stein}. 
\medskip\\
\textbf{Acknowledgements.} The author would like to thank Anar Akhmedov, Mikio Furuta, Kazunori Kikuchi, Hirofumi Sasahira, Motoo Tange,  Masaaki Ue and Yuichi Yamada for their helpful comments. 
\section{Log transform, knot surgery and the adjunction inequality}\label{sec:log}
In this section, we briefly recall log transform, knot surgery and their effects on Seiberg-Witten invariants. Furthermore, we recall the adjunction inequality. For details, see \cite{GS}, \cite{FS4}, \cite{S}, \cite{Sco}. 

 Throughout this paper, for an embedded oriented connected closed surface $\Sigma$ in a $4$-manifold $X$, we denote by $[\Sigma]$ the class of $H_2(X;\mathbb{Z})$ represented by $\Sigma$. The symbol $\nu (\cdot)$ and $PD(\cdot)$ denote the regular neighborhood and the Poincar\'{e} dual class, respectively. 

\begin{figure}[ht!]
\begin{center}
\includegraphics[width=4.5in]{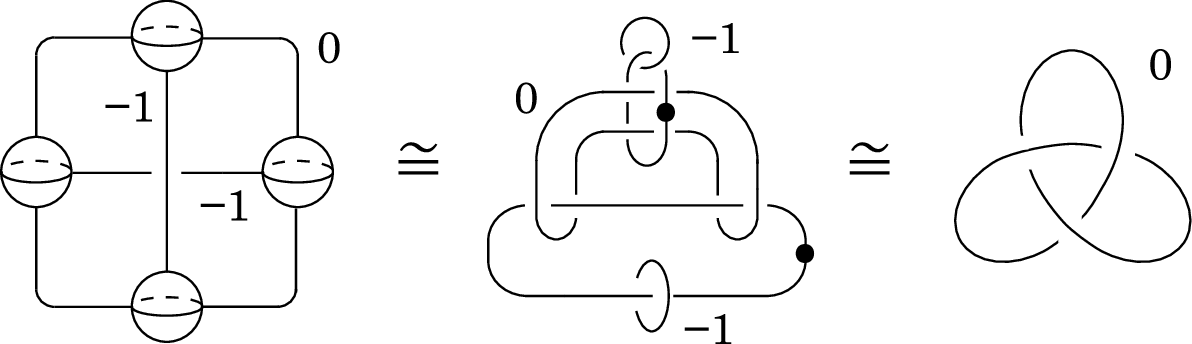}
\caption{Diagrams of a cusp neighborhood $C$}
\label{fig:nuclei_fig1}
\end{center}
\end{figure}

\begin{definition}$(1)$ The smooth $4$-manifold $C$ given by the left handlebody in Figure~\ref{fig:nuclei_fig1} is called a \textit{cusp neighborhood}. Note that this $4$-manifold is diffeomorphic to the middle and the right handlebodies. It is known that $C$ admits a torus fibration structure over $D^2$ with only one singular fiber, i.e., a cusp fiber. The obvious $T^2\times D^2$ in the left diagram describes the tubular neighborhood of the regular fiber. \\
$(2)$ A smoothly embedded torus $T$ in a $4$-manifold $X$ is called a \textit{$c$-embedded torus} if $X$ smoothly contains a cusp neighborhood $C$, and $T$ is a regular fiber of $C$. Note that the self-intersection number of any $c$-embedded torus is zero. 
\end{definition}
\subsection{Log transform}Recall that self-diffeomorphisms of $T^3=S^1\times S^1\times S^1$ are classified, up to isotopy, by $\textnormal{Aut}_\mathbb{Z}(H_1(T^3;\mathbb{Z}))$.  Let $\varphi_p$ $(p\geq 0)$ be the self-diffeomorphism of $T^2\times S^1$ induced by the automorphism 
\begin{equation*}
\left(
\begin{array}{ccc}
1 &0 &0  \\
0 &0 &1  \\
0 &-1 &p 
\end{array}
\right)
\end{equation*}
of $H_1(S^1;\mathbb{Z})\oplus H_1(S^1;\mathbb{Z})\oplus H_1(S^1;\mathbb{Z})$ with the obvious basis. 
\begin{definition}Let $X$ be a smooth $4$-manifold which contains a torus $T$ with the self-intersection number zero, and let $\nu (T)$ denote the tubular neighborhood of $T$. Let $X_{(p)}$ be the smooth $4$-manifold obtained from $X$ by removing $\nu (T)$ and gluing $T^2\times D^2$ via the diffeomorphism $\varphi_p: T^2\times \partial D^2\to \partial \nu(T)$, where we identify $\partial \nu(T)$ with $T^2\times \partial D^2$. This $X_{(p)}$ is called the \textit{$p$-log transform} of $X$ along $T$. 
We denote the new embedded torus $T^2\times {0}$ in $X_{(p)}$ by $T_{p}$, which is called the \textit{multiple fiber}. 
\end{definition}
\begin{remark}\label{rem:multiple fiber}$(1)$ Note that $X_{(p)}$ contains a parallel copy of the original torus $T$. This torus satisfies $[T]=p[T_p]$ in $H_2(X_{(p)};\mathbb{Z})$. This can be seen as follows. Let $\lambda_1, \lambda_2, \lambda_3$ be the obvious basis of $H_1(S^1\times S^1;\mathbb{Z})\oplus H_1(\partial D^2;\mathbb{Z})$. 
The tori $\varphi_p(T^2\times \{pt.\})$ and $T\times \{pt.\}$ in $\partial \nu(T)\subset X_{(p)}$ are parallel copies of $T_p$ and $T$, respectively. By the definition of $\varphi_p$, the obvious basis $\mu_1, \mu_2, \mu_3$ of $H_1(\partial\nu(T))=H_1(T;\mathbb{Z})\oplus H_1(\partial D^2;\mathbb{Z})$ satisfies $\mu_1={\varphi_p}_*(\lambda_1)$ and $\mu_2=p{\varphi_p}_*(\lambda_2)-{\varphi_p}_*(\lambda_3)$. Since $[\varphi_p(S^1\times\{pt.\}\times S^1)]=0$ in $H_2(X_{(p)};\mathbb{Z})$, the K\"{u}nuth formula easily gives $[T]=p[T_p]$ in $H_2(X_{(p)};\mathbb{Z})$ as desired. 
\\
$(2)$ The above definition of log transform is a special case of the general log transform. When the torus is $c$-embedded, this definition is known to be consistent with the general one. See~\cite{G0}, \cite{GS} for details. 
\end{remark}
By the following result, we may identify $X_{(1)}$ with $X$ if the torus $T$ is $c$-embedded. 
\begin{theorem}[Gompf~\cite{G0}] Let $X$ be a connected closed oriented smooth $4$-manifold which contains a $c$-embedded torus $T$. Then the log transform $X_{(1)}$ along $T$ is diffeomorphic to $X$. 
\end{theorem}
The Seiberg-Witten invariant of $X_{(p)}$ is given as follows. 
\begin{theorem}[Fintushel-Stern~\cite{FS1}, Morgan-Mrowka-Szab\'{o}~\cite{MMS}]\label{Thm:SW:log}Let $X$ be a connected closed oriented smooth $4$-manifold with $b_2^+>1$ which contains a $c$-embedded torus $T$. Let $X_{(p)}$ and $T_p$ $(p\geq 1)$ be the $p$-log transform of $X$ along $T$ and the multiple fiber of $X_{(p)}$, respectively. If $[T]$ is non-torsion in $H_2(X;\mathbb{Z})$, then 
\begin{equation*}
\mathcal{SW}_{X_{(p)}}=\mathcal{SW}_X\cdot (t^{-(p-1)}+t^{-(p-3)}+\dots+t^{p-1}), 
\end{equation*}
where $t=\exp(PD([T_{p}]))$. 
\end{theorem}

\subsection{Knot surgery}
\begin{definition}\label{def:knot surgery}Let $X$ be a compact oriented smooth $4$-manifold which contains a torus $T$ with the self-intersection number $0$, and let $K$ be a knot in $S^3$. We call the smooth $4$-manifold 
\begin{equation*}
X_{K}:=(X-\text{int}\, \nu (T))\cup_{\varphi} ((S^3-\text{int}\, \nu (K))\times S^1)
\end{equation*} 
a \textit{knot surgery} of $X$ along $T$ with $K$. Here ${\varphi}: \partial \nu (T)\to \partial \nu (K)\times S^1$ is any gluing map which sends $\{pt.\}\times \partial D^2$ $(\subset T^2\times \partial{D^2}\cong \partial\nu (T))$ to $\ell\times \{pt.\}(\subset \partial \nu (K)\times S^1)$, where $\ell$ denotes a longitude of $K$. Note that the diffeomorphism type of $X_K$ may depend on the choice of ${\varphi}$.  We denote by $T_K$ the torus $m\times S^1\subset (S^3-\nu(K))\times S^1$ in $X_K$, where $m$ denotes a meridian of $K$. Note $[T_K]\cdot [T_K]=0$. 
\end{definition}
\begin{remark}\label{rem:knot:T_K}
Note that $X_K$ contains a parallel copy of the original torus $T$. Similarly to Remark~\ref{rem:multiple fiber}.(1), we can see that $[T]=[T_K]$ in $H_2(X_K;\mathbb{Z})$. 
\end{remark}
The Seiberg-Witten invariant of $X_{K}$ is given as follows. 
\begin{theorem}[Fintushel-Stern~\cite{FS2}, see also \cite{FS5}]\label{Thm:SW:knot}
Let $X$ be a connected closed oriented smooth $4$-manifold with $b_2^+>1$ which contains a $c$-embedded torus $T$. If $[T]$ is non-torsion in $H_2(X;\mathbb{Z})$, then, for any knot $K$ in $S^3$, the knot surgery $X_K$ satisfies  
\begin{equation*}
\mathcal{SW}_{X_K}=\mathcal{SW}_X\cdot \Delta_K(t), 
\end{equation*}
where $\Delta_K$ is the symmetrized Alexander polynomial of K, and $t=\exp(PD(2[T]))$. 
\end{theorem}
\subsection{The adjunction inequality}In this paper, we use the following adjunction inequality to detect smooth structures. According to Theorem~\ref{thm:simple type}, we can apply the adjunction inequality for many closed $4$-manifolds with the non-vanishing Seiberg-Witten invariants. 
\begin{theorem}[Kronheimer-Mrowka~\cite{KM1}, Ozsv\'{a}th-Szab\'{o}~\cite{OzSz}. cf.~\cite{OS1}]Let $X$ be a connected closed oriented smooth $4$-manifold with $b_2^+ >1$,
 and let $K\in H^2(X;\mathbb{Z})$ be a Seiberg-Witten basic class of $X$. If a smoothly embedded connected closed oriented surface $\Sigma\subset X$ of genus $g\geq 0$ satisfies $[\Sigma]\cdot [\Sigma]\geq 0$, and $[\Sigma]$ is non-torsion, then the following inequality holds: 
\begin{equation*}
[\Sigma]\cdot [\Sigma]+\lvert \langle K, [\Sigma] \rangle\rvert\leq 2g-2.
\end{equation*}
Furthermore, if $X$ is of Seiberg-Witten simple type and $g\geq 1$, then the same inequality also holds in the $[\Sigma]\cdot [\Sigma]<0$ case. 
\end{theorem}
\begin{theorem}[Morgan-Mrowka-Szabo~\cite{MMS}]\label{thm:simple type}
Let $X$ be a connected closed oriented smooth $4$-manifold with $b_2^+ >1$. If $X$ contains a smoothly embedded torus $T$, and $[T]$ is non-torsion, then $X$ is of Seiberg-Witten simple type. 
\end{theorem}

\section{Generalized nucleus}\label{sec:nucleus}In this section, we give a new generalization of Gompf nuclei~\cite{G0} and discuss their properties. For related discussions, see also~\cite{FS0}. 
We begin with the following lemma. 
\begin{lemma}\label{nucleus:lem:divisor}Let $X$ be a compact connected oriented smooth $4$-manifold $($possibly with boundary$)$ which contains a torus $T$. Suppose that $H_2(X;\mathbb{Z})$ has no torsion and that $H_1(X-\textnormal{int}\, \nu (T);\mathbb{Z})\cong \mathbb{Z}/d\mathbb{Z}$ for some positive integer $d$. Then there exists a primitive $($i.e.\ indivisible$)$ class $\widehat{T}$ of $H_2(X;\mathbb{Z})$ such that $[T]=d\cdot \widehat{T}$.
\end{lemma}
\begin{proof}Using the homology exact sequence for the pair $(X, \nu (T))$, the Alexander duality $H^{4-i}(X-\text{int}\, \nu (T);\mathbb{Z})\cong H_i(X,\nu (T);\mathbb{Z})$, and the universal coefficient theorem, we easily get the following exact sequence:
\begin{equation*}
0\to H_2(\nu (T);\mathbb{Z})\to H_2(X;\mathbb{Z})\to \textnormal{Hom}(H_2(X-\textnormal{int}\, \nu(T);\mathbb{Z}), \mathbb{Z})\oplus (\mathbb{Z}/d\mathbb{Z}) \to \mathbb{Z}\oplus \mathbb{Z}.
\end{equation*}
We can now easily check the claim. 
\end{proof}
We here give a generalization of Gompf nucleus~\cite{G0}. Our generalization is different from Fuller's one~\cite{Fu}, since our purpose is to produce exotic smooth structures. 
\begin{definition}\label{def:nucleus}Let $N$ be a connected compact oriented smooth $4$-manifold with a connected non-empty boundary. Suppose that $N$ smoothly contains a torus $T$. We call $(N,T)$ (or $N$ itself) \textit{a nucleus}, if the following conditions hold. 
\begin{itemize}
\item [(i)] $\pi_1(N)\cong 1$ and $H_2(N;\mathbb{Z})\cong \mathbb{Z}\oplus \mathbb{Z}$.
 \item [(ii)] The intersection form of $N$ is unimodular. In particular, the boundary $\partial N$ is a connected homology $3$-sphere (cf.\ \cite{GS}). 
 \item [(iii)] $T$ is $c$-embedded in $N$. In particular, $[T]\cdot [T]=0$. 
 \item [(iv)] $\pi_1(N-\text{int}\,\nu (T))\cong \mathbb{Z}/d_T\mathbb{Z}$ for some positive integer $d_T$. We call $d_T$ the \textit{divisor} of $T$. 
 \item [(v)] The inclusion induced homomorphism $\pi_1(\partial \nu (T))\to \pi_1(N-\text{int}\, \nu (T))$ is surjective. 
\end{itemize}
\end{definition}
\begin{remark}
In the above definition of nuclei, we added the condition (v) to guarantee the simple connectivity of the log transform $N_{(p)}$ for some $p$'s (see Lemma~\ref{sec:nucleu:lem:property}). However, as long as $N_{(p)}$ is simply connected (for some $p$'s),  all arguments in this paper work even when $N$ does not satisfy the condition (v).
\end{remark}
There are many examples of our generalized nuclei. Actually, it is easy to construct examples using handlebody pictures. 
\begin{example}$(1)$ Let $X$ be a $4$-manifold obtained from the cusp neighborhood $C$ by attaching a single $2$-handle along a knot which links with the 0-framed trefoil knot in the right diagram of Figure~\ref{fig:nuclei_fig1} geometrically once. Let $T$ denote the regular fiber of the cusp neighborhood $C$. Then $(X, T)$ is clearly a nucleus with $d_T=1$. It easily follows from the Van Kampen's Theorem that $(X_{(p)}, T)$ is also a nucleus with $d_T=p$. Note that $X_{(p)}$ still contains a smaller copy of $C$. \smallskip\\
$(2)$ Let $G(n)$ $(n\geq 1)$ be the $4$-manifold in Figure~\ref{fig:nuclei_fig2}. This $G(n)$ is the Gompf nucleus introduced by Gompf~\cite{G0}. Note that the intersection form of $G(n)$ is odd (resp.\ even) if $n$ is odd (resp.\ even). As noted in the above, $G(n)$ and its $p$-log transform ${G(n)}_{(p)}$ $(p\geq 1)$ are nuclei in the sense of Definition~\ref{def:nucleus}. 
\begin{figure}[ht!]
\begin{center}
\includegraphics[width=1.0in]{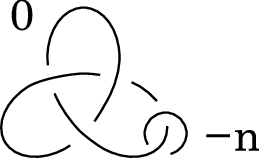}
\caption{Gompf nucleus $G(n)$ $(n\geq 1)$}
\label{fig:nuclei_fig2}
\end{center}
\end{figure}
\end{example}

\begin{figure}[ht!]
\begin{center}
\includegraphics[width=4.3in]{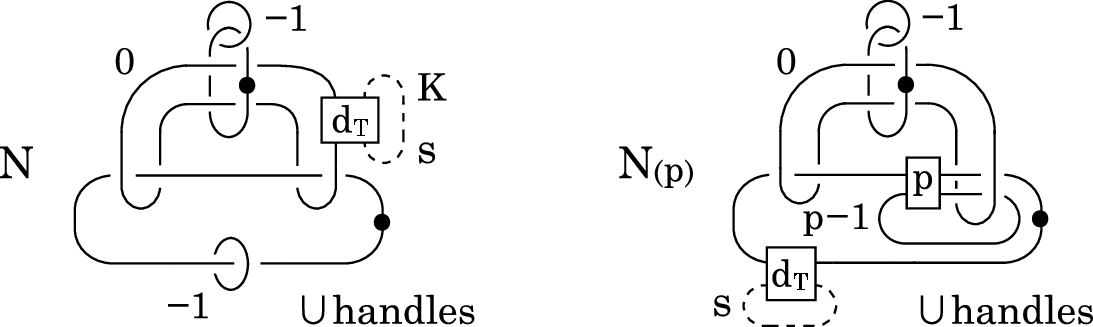}
\caption{handlebodies of $N$ and $N_{(p)}$}
\label{fig:nuclei_fig3}
\end{center}
\end{figure}

We here discuss basic properties of our nuclei.  
\begin{lemma}\label{nuclei:lem:basic:log}Let $(N,T)$ and $d_T$ be a nucleus and the divisor of $T$, respectively. Then there exists a basis $\widehat{T}, S$ of $H_2(N;\mathbb{Z})$ such that $[T]=d_T\cdot \widehat{T}$ and $S\cdot \widehat{T}=1$. Consequently, the intersection form of $N$ is indefinite. 
\end{lemma}
\begin{proof}Lemma~\ref{nucleus:lem:divisor} gives a primitive class $\widehat{T}$ with $[T]=d_T\cdot \widehat{T}$. Since the intersection form of $N$ is unimodular, the fact $\widehat{T}\cdot \widehat{T}=0$ implies the existence of the required class $S$. 
\end{proof}
\begin{lemma}\label{sec:nucleu:lem:property}Let $(N,T)$ and $d_T$ be a nucleus and the divisor of $T$, respectively. Assume that a positive integer $p$ is relatively prime to $d_T$. Let $N_{(p)}$ be the $p$-log transform of $N$ along $T$. Then the following hold.\\
$(1)$ $\pi_1(N_{(p)})\cong 1$ and $H_2(N_{(p)};\mathbb{Z})\cong \mathbb{Z}\oplus \mathbb{Z}$.
\\
$(2)$ The intersection form of $N_{(p)}$ is unimodular and indefinite. The intersection form of $N_{(p)}$ is even if and only if the following two conditions are satisfied. 
\begin{enumerate}
 \item [(i)] The intersection form of $N$ is even. 
 \item [(ii)] $p$ is odd, or $d_T$ is even. 
\end{enumerate}
$(3)$ If the intersection form of $N$ and $N_{(p)}$ are isomorphic, then any diffeomorphism between the boundaries $\partial N$ and $\partial N_{(p)}$ extends to a homeomorphism between $N$ and $N_{(p)}$. 
\end{lemma}
\begin{proof}$(1)$ The Van Kampen's theorem gives $\pi_1(N_{(p)})\cong 1$. We can easily check that $H_2(N_{(p)}; \mathbb{Z})$ has no torsion and that $H_3(N_{(p)};\mathbb{Z})=0$, using the universal coefficient theorem and the Poincar\'{e} duality. Since the Euler characteristics of $N_{(p)}$ and $N$ are equal, $H_2(N_{(p)};\mathbb{Z})\cong \mathbb{Z}\oplus \mathbb{Z}$. 

$(2)$ Since $N$ is simply connected, and $\partial N_{(p)}(=\partial N)$ is a homology $3$-sphere, the intersection form of $N_{(p)}$ is unimodular (cf.\ \cite{GS}). Lemma~\ref{nucleus:lem:divisor} thus implies that $N_{(p)}$ has an indefinite form. 

We next discuss the parity of the intersection form. Let $S$ denote the class of $H_2(N;\mathbb{Z})$ given in Lemma~\ref{nuclei:lem:basic:log} and put $s=S\cdot S$. We here consider a handle decomposition of $N$. $N$ is obtained from the cusp neighborhood in the middle diagram of Figure~\ref{fig:nuclei_fig1} by attaching handles. 
Recall that (part of) 2-handles of $N$ naturally give a generating set of $H_2(N;\mathbb{Z})$ after appropriately sliding 2-handles.  Therefore, if necessary by creating a canceling 2-handle/3-handle pair and sliding handles except those of the cusp neighborhood $C$, the class $S$ is represented by a $2$-handle (say $K$) of the decomposition of $N$.  In particular, (the attaching circle of) $K$ does not algebraically link with any $1$-handles of the decomposition of $N$. Furthermore, we may assume that $K$ does not algebraically link with the two $-1$-framed knots of $C$, if necessary, by further sliding $K$ over $1$-handles of $C$. It follows that the Seifert framing of $K$ is $s$ and that $K$ links with the 0-framed knot in the left picture of Figure~\ref{fig:nuclei_fig3} algebraically $d_T$ times. Hence the $2$-handle $K$ is given as shown in the left picture of Figure~\ref{fig:nuclei_fig3}, though the actual $K$ may geometrically link with other handles in the picture. 

Applying the $p$-log transform procedure in Figure~17 of~\cite{AY1} to the above handlebody of $N$, we get the right handlebody of $N_{(p)}$ in Figure~\ref{fig:nuclei_fig3}. Note that the obvious $T^2\times D^2$ describes $\nu(T_p)$. The bottom $s$-framed circle (say $L_1$) may geometrically link with other handles in the picture, though its algebraic linking numbers with other handles are correct. Let $L_2$ denote the middle $(p-1)$-framed knot in the right diagram of Figure~\ref{fig:nuclei_fig3}. By introducing a canceling $2$-handle/$3$-handle pair and sliding handles, we get a $2$-handle $L$ $(=pL_1+d_TL_2)$ which satisfies the following three conditions: (i) $L$ does not algebraically link with any $1$-handles; (ii) $L$ links with the middle $0$-framed knot algebraically $d_T$ times; (iii) The Seifert framing of $L$ is $s':=p^2s+(d_T)^2(p-1)$. This together with Lemma~\ref{nucleus:lem:divisor} implies the existence of a basis $S', \widehat{T}_p$ of $H_2(N_{(p)};\mathbb{Z})$ which satisfies the following:
\begin{equation*}
[T_p]=d_T\cdot \widehat{T}_p,\quad S'\cdot S'=s',\quad S'\cdot \widehat{T}_p=1. 
\end{equation*}
Thus the intersection form of $N_{(p)}$ is even if and only if $s'=p^2s+(d_T)^2(p-1)$ is even. We can now easily check the claim (2). 

$(3)$ This clearly follows from Boyer's theorem~\cite{B} (cf.\ Corollary~2.3 of \cite{P2}). 
\end{proof}
\begin{lemma}\label{lem:basic:nuclei:knot}Let $(N,T)$ be a nucleus. Assume that the divisor $d_T$ satisfies $d_T=1$. 
Let $K$ be a knot in $S^3$, and let $N_K$ be the knot surgery of $N$ along $T$. Then the following hold.\\
$(1)$ The fundamental group, the homology groups, and the intersection form of $N_{K}$ are isomorphic to those of $N$. Furthermore, the class $[T]$ is primitive in $H_2(N_K;\mathbb{Z})$. 
\\
$(2)$ Any diffeomorphism between the boundaries $\partial N$ and $\partial N_K$ extends to a homeomorphism between $N$ and $N_K$. 
\end{lemma}
\begin{proof}$(1)$ A sketch of a proof was given by Fintushel-Stern~\cite{FS2}. We here explain more detail for the reader's convenience. Since $\pi_1(X-\nu(T))\cong 1$, $N_K$ is simply connected (\cite{FS2}). Since $\partial N_K$ is a connected homology $3$-sphere, it is easy to see that $H_2(N_K;\mathbb{Z})$ has no torsion and that $H_3(N_K;\mathbb{Z})=0$ by the Poincar\'{e} duality and the universal coefficient theorem. These fact give $H_2(N_K;\mathbb{Z})\cong H_2(N;\mathbb{Z})$ because the Euler characteristics of $N_K$ and $N$ are the same. Hence the homology groups of $N_K$ and $N$ are isomorphic. 

The condition $\pi_1(X-\nu(T))\cong 1$ gives an embedded closed connected oriented surface $S$ in $N$ which transversely intersects with the torus $T$ in exactly one point. 
(The reason is as follows. It is easy to construct a continuous map $S^2\to N$ such that its image of a disk $(\subset S^2)$ is a disk $D:=\{pt.\}\times D^2$ in $T\times D^2\cong \nu(T)$. Approximating this map by an immersion and resolving its singular points (cf.\ Remark~2.1.2 of \cite{GS}), we obtain the desired surface.) 
Let $T'\subset N-\nu(T)$ be a parallel copy of $T$. We may assume that $T'\subset N$ intersects with $S\subset N$ in exactly one point and that the point is outside of the disk $D$. 
Let $\Sigma$ be a Seifert surface of a longitude of $K$. We denote by $S_K$ the closed connected oriented surface in $X_K$ obtained by gluing $S-\textnormal{int}\, D$ with $\Sigma\times \{pt.\}\subset (S^3-\textnormal{int}\, \nu(K))\times S^1$ along $\partial D=\partial \Sigma \times \{pt.\}$. 
Since we can move $\Sigma$ along the trivial $S^1$ direction in $(S^3-\textnormal{int}\, \nu(K))\times S^1$, the surface $S_K$ satisfies $S_K\cdot S_K=S\cdot S$. By the construction, $S_K\cdot T'=1$. This shows that the intersection forms of $N$ and $N_K$ are isomorphic and that the class $[T]$ is primitive in $H_2(N_K;\mathbb{Z})$. 

The claim (2) follows from Boyer's theorem. 
\end{proof}
\section{Relative genus function}\label{sec:minimal genus}In this section, we introduce a relative genus function of a $4$-manifold. 

Here recall that, for any oriented smooth $4$-manifold $X$, every class $\alpha\in H_2(X;\mathbb{Z})$ can be represented by a smoothly embedded closed connected oriented surface in $X$ (cf.\ Exercise~4.5.12.(b) of \cite{GS}). Note that this holds even when $X$ has (possibly disconnected) boundary, or is non-compact. We also use the following convention. For a class $\alpha\in H_2(X;\mathbb{Z})/\textnormal{Tor}$ and a smoothly embedded closed connected oriented surface $\Sigma$ in $X$, we say ``$\Sigma$ represents $\alpha$'', if $[\Sigma]\in H_2(X;\mathbb{Z})$ represents $\alpha\in H_2(X;\mathbb{Z})/\textnormal{Tor}$. 

The (minimal) genus function $g_X:H_2(X;\mathbb{Z})\to \mathbb{Z}$ of $X$ is defined as follows: 
\begin{equation*}
g_X(\alpha) = \min\left\{g\in \mathbb{Z}\, \left| \begin{tabular}{l}
\text{$\alpha$ is represented by a smoothly embedded}\\\text{closed connected oriented surface of genus $g$ in $X$}.
\end{tabular}\right.\right\}.
\end{equation*}
Similarly, one can also define the genus function $g_X:H_2(X;\mathbb{Z})/\textnormal{Tor}\to \mathbb{Z}$ of $X$. We use the same symbol $g_X$ for these two functions, if it is clear from the context.

In~\cite{A2}, Akbulut successfully used the genus function to detect smooth structures of two $4$-manifolds with $b_2=1$. However, it is generally difficult to see differences of the functions in the $b_2\geq 2$ case, since we have to consider the automorphisms of the second homology groups induced by self-diffeomorphisms. 

To avoid this issue, we introduce relative versions of the genus function, by generalizing techniques of Akbulut and the author in \cite{AY5}. 
We use the following notation. Let $\mathbb{Z}_{\geq 0}$ denote the set of non-negative integers. For an index set $\Lambda$, let $\mathbb{Z}^\Lambda$ and $(\mathbb{Z}_{\geq 0})^\Lambda$ denote the direct products $\prod _{\Lambda} \mathbb{Z}$ and $\prod _{\Lambda} (\mathbb{Z}_{\geq 0})$, respectively. Let $\textnormal{Sym}_\Lambda(\mathbb{Z})$ denote the set of $\Lambda\times \Lambda$ integer symmetric matrices. Namely, we set
\begin{align*}
\mathbb{Z}^\Lambda &= \{(d_\lambda)_{\lambda\in \Lambda}\mid \text{$d_\lambda \in \mathbb{Z}$ for all $\lambda\in \Lambda$}\},
\\
(\mathbb{Z}_{\geq 0})^\Lambda &= \{(g_\lambda)_{\lambda\in \Lambda}\mid \text{$g_\lambda \in \mathbb{Z}_{\geq 0}$ for all $\lambda\in \Lambda$}\}, 
\\
\textnormal{Sym}_\Lambda(\mathbb{Z}) &= \{(Q_{\lambda,\mu})_{(\lambda,\mu)\in \Lambda\times \Lambda}\mid \text{$Q_{\lambda, \mu}\in \mathbb{Z}$ and $Q_{\lambda, \mu}=Q_{\mu,\lambda}$ for all $(\lambda,\mu)\in \Lambda\times \Lambda$}\}.
\end{align*}

Let $X$ be an arbitrary connected oriented smooth $4$-manifold. $X$ is allowed to be closed or to have (possibly disconnected) boundary or to be non-compact. 
\begin{definition}\label{def:relative genus}
Suppose that $H_2(X;\mathbb{Z})/\textnormal{Tor}$ is a free $\mathbb{Z}$-module. Let fix an index set $\Lambda$ which has a bijection to a basis of $H_2(X;\mathbb{Z})/\textnormal{Tor}$. Note that any two bases of a (possibly infinitely generated) free $\mathbb{Z}$-module have the same cardinality. Fix an element $\lambda_0$ of $\Lambda$. 
\smallskip\\
(1) Let $\mathbf{v}=\{v_\lambda\mid  \lambda \in \Lambda \}$ be a basis of $H_2(X;\mathbb{Z})/\textnormal{Tor}$. For $Q\in \textnormal{Sym}_\Lambda(\mathbb{Z})$, $\mathbf{d}\in \mathbb{Z}^\Lambda$ and $\mathbf{g}\in (\mathbb{Z}_{\geq 0})^{\Lambda-\{\lambda_0\}}$, we define $G_{X, Q, \mathbf{d}, \mathbf{g}}(\mathbf{v})\in \mathbb{Z}\cup \{\infty\}$ as follows. Here we put $Q=(Q_{\lambda,\mu})_{(\lambda,\mu)\in \Lambda\times \Lambda}$, $\mathbf{d}=(d_\lambda)_{\lambda\in \Lambda}$ and $\mathbf{g}=(g_\lambda)_{\lambda\in \Lambda-\{\lambda_0\}}$.  
\begin{equation*}
G_{X, Q, \mathbf{d}, \mathbf{g}}(\mathbf{v})=
\begin{cases}
g_X(d_{\lambda_0}v_{\lambda_0}),&\text{if $\mathbf{v}$ satisfies the following conditions (i) and (ii)}.
\\
\infty,&\text{otherwise}. 
\end{cases}
\end{equation*}

\begin{enumerate}
 \item [(i)] The basis $\mathbf{v}$ represents $Q$, namely, $v_{\lambda}\cdot v_{\mu}=Q_{\lambda,\mu}$ for each $\lambda, \mu\in \Lambda$. 
 \smallskip
 \item [(ii)] For each $\lambda\in \Lambda-\{\lambda_0\}$, $d_\lambda v_\lambda\in H_2(X;\mathbb{Z})/\textnormal{Tor}$ satisfies $g_X(d_\lambda v_\lambda)\leq g_\lambda$. 
\end{enumerate}
(2)
For $Q\in \textnormal{Sym}_\Lambda(\mathbb{Z})$, $\mathbf{d}\in \mathbb{Z}^\Lambda$ and $\mathbf{g}\in (\mathbb{Z}_{\geq 0})^{\Lambda-\{\lambda_0\}}$, we define $G_X(Q, \mathbf{d}, \mathbf{g})\in \mathbb{Z}\cup \{\infty\}$ as follows. 
\begin{equation*}
G_X(Q, \mathbf{d}, \mathbf{g})=\min\{G_{X, Q, \mathbf{d}, \mathbf{g}}(\mathbf{v})\mid \text{$\mathbf{v}$ is a basis of $H_2(X;\mathbb{Z})/\textnormal{Tor}$.}\}.
\end{equation*}
We call this function
\begin{equation*}
G_X: \textnormal{Sym}_\Lambda(\mathbb{Z})\times {\mathbb{Z}}^\Lambda\times (\mathbb{Z}_{\geq 0})^{\Lambda-\{\lambda_0\}}\to \mathbb{Z}\cup \{\infty\}
\end{equation*}
the \textit{relative genus function} of $X$.
\end{definition}
\begin{remark}$(1)$ 
This function is clearly an invariant of smooth structures. Namely, if $X$ and $Y$ are orientation preserving diffeomorphic, then $G_{X}=G_{Y}$ for any fixed pair $(\Lambda, \lambda_0)$. If $X$ and $Y$ are orientation reversing diffeomorphic, then $G_X(Q, \mathbf{d}, \mathbf{g})=G_Y(-Q, \mathbf{d}, \bf{g})$ for each $(Q, \mathbf{d}, \mathbf{g})\in \textnormal{Sym}_\Lambda(\mathbb{Z})\times {\mathbb{Z}}^\Lambda\times (\mathbb{Z}_{\geq 0})^{\Lambda-\{\lambda_0\}}$.
\smallskip\\
$(2)$ The following variant $G^\pm$ is useful to detect smooth structures without fixing orientations. Let\begin{equation*}
G^{\pm}_X: \textnormal{Sym}_\Lambda(\mathbb{Z})\times {\mathbb{Z}}^\Lambda\times (\mathbb{Z}_{\geq 0})^{\Lambda-\{\lambda_0\}}\to \mathbb{Z}\cup \{\infty\}
\end{equation*} 
be the function defined by $G^\pm_X(Q, \mathbf{d}, \mathbf{g})=\min \{G_X(Q, \mathbf{d}, \mathbf{g}),\, G_X(-Q, \mathbf{d}, \mathbf{g})\}$. Clearly, if $X$ and $Y$ are diffeomorphic (not necessarily orientation preserving), then $G^{\pm}_X=G^{\pm}_Y$.
\end{remark}

We next slightly generalize the relative genus function, since we need to deal with the case where $H_2(X;\mathbb{Z})/\textnormal{Tor}$ is not a free module. While the following function $G_{X, F_0, F_1}$ is stronger and enough for our purpose, we mainly use $G_{X}$ in this paper. This is just because we can decrease the number of symbols for discussions. 
\begin{definition}\label{def:relative genus:minor}
Let $F_0$ and $F_1$ be a free $\mathbb{Z}$-module and a (possibly non-free) $\mathbb{Z}$-module, respectively. Let $\Lambda_0$ be an index set which has a bijection to a basis of $F_0$. 
 Put $\Lambda_1=F_1$ and regard $\Lambda_1$ as an index set. Denote by $\Lambda$ the disjoint union $\Lambda_0 \coprod \Lambda_1$. Fix an element $\lambda_0\in \Lambda_0\subset \Lambda$. 
\smallskip\\
$(1)$ Let $\mathbf{v}=\{v_\lambda\mid  \lambda \in \Lambda \}$ be a subset of $H_2(X;\mathbb{Z})$. 
For $Q\in \textnormal{Sym}_{\Lambda}(\mathbb{Z})$, $\mathbf{d}\in \mathbb{Z}^\Lambda$ and $\mathbf{g}\in (\mathbb{Z}_{\geq 0})^{\Lambda-\{\lambda_0\}}$, we define $G_{X, F_0, F_1,Q, \mathbf{d}, \mathbf{g}}(\mathbf{v})\in \mathbb{Z}\cup \{\infty\}$ as follows. Here we put $Q=(Q_{\lambda,\mu})_{(\lambda,\mu)\in \Lambda\times \Lambda}$, $\mathbf{d}=(d_\lambda)_{\lambda\in \Lambda}$ and $\mathbf{g}=(g_\lambda)_{\lambda\in \Lambda-\{\lambda_0\}}$.
\begin{equation*}
G_{X, F_0, F_1,Q, \mathbf{d}, \mathbf{g}}(\mathbf{v})=
\begin{cases}
g_X(d_{\lambda_0}v_{\lambda_0}),&\text{if $\mathbf{v}$ satisfies the following conditions (i)--(iii).}
\\
\infty,&\text{otherwise}. 
\end{cases}
\end{equation*}
\begin{enumerate}
  \item [(i)] The subset $\mathbf{v}$ represents $Q$, namely, $v_{\lambda}\cdot v_{\mu}=Q_{\lambda,\mu}$ for each $\lambda, \mu\in \Lambda$. 
  \item [(ii)] For each $\lambda\in \Lambda-\{\lambda_0\}$, $d_\lambda v_\lambda\in H_2(X;\mathbb{Z})$ satisfies $g_X(d_\lambda v_\lambda)\leq g_\lambda$. 
 \item [(iii)] There exists a direct sum decomposition $H_2(X;\mathbb{Z})=F_0'\oplus F_1'$ satisfying the following conditions (a) and (b).
 \begin{enumerate}
 \item [(a)] The submodules $F_0'$ and $F_1'$ are isomorphic to $F_0$ and $F_1$, respectively.
 \item [(b)] The subset $\mathbf{v}^{(0)}:=\{v_\lambda\mid  \lambda \in \Lambda_0 \}$ of $\mathbf{v}$ is a basis of $F_0'$, and the subset $\mathbf{v}^{(1)}:=\{v_\lambda\mid  \lambda \in \Lambda_1 \}$ of $\mathbf{v}$ satisfies $\mathbf{v}^{(1)}=F_1'$. 
\end{enumerate}
\end{enumerate}
(2) For $Q\in \textnormal{Sym}_{\Lambda}(\mathbb{Z})$, $\mathbf{d}\in \mathbb{Z}^\Lambda$ and $\mathbf{g}\in (\mathbb{Z}_{\geq 0})^{\Lambda-\{\lambda_0\}}$, we define $G_{X,F_0,F_1}(Q, \mathbf{d}, \mathbf{g})\in \mathbb{Z}\cup \{\infty\}$ as follows.
\begin{equation*}
G_{X,F_0,F_1}(Q, \mathbf{d}, \mathbf{g})=\min\{G_{X, F_0, F_1,Q, \mathbf{d}, \mathbf{g}}(\mathbf{v})\mid \text{$\mathbf{v}$ is a subset of $H_2(X;\mathbb{Z})$.}\}.
\end{equation*}
We call this function
\begin{equation*}
G_{X, F_0, F_1}: \textnormal{Sym}_\Lambda(\mathbb{Z})\times {\mathbb{Z}}^\Lambda\times (\mathbb{Z}_{\geq 0})^{\Lambda-\{\lambda_0\}}\to \mathbb{Z}\cup \{\infty\}
\end{equation*}
 the \textit{relative genus function} of $X$ with respect to $F_0, F_1$.
\end{definition}
\begin{remark}
Similarly to $G_X$, this function is an invariant of smooth structures. Namely, if $X$ and $Y$ are orientation preserving diffeomorphic, then $G_{X, F_0, F_1}=G_{Y, F_0, F_1}$ for each modules $F_0$ and $F_1$. If $X$ and $Y$ are orientation reversing diffeomorphic, then $G_{X, F_0, F_1}(Q, \mathbf{d}, \mathbf{g})=G_{Y, F_0, F_1}(-Q, \mathbf{d}, \mathbf{g})$ for each $(Q, \mathbf{d}, \mathbf{g})\in \textnormal{Sym}_\Lambda(\mathbb{Z})\times {\mathbb{Z}}^\Lambda\times (\mathbb{Z}_{\geq 0})^{\Lambda-\{\lambda_0\}}$. Similarly to $G^{\pm}_X$, we can also define $G^\pm_{X, F_0, F_1}$. 
\end{remark}

\section{Producing exotic smooth structures}\label{sec:nuclei and exotic}

To state the results of this paper simply, we use the following convention. For a smooth $4$-manifold $X$, we say ``$X$ admits infinitely many distinct (exotic) smooth structures'', if there exist infinitely many smooth $4$-manifolds which are all homeomorphic to $X$ but mutually non-diffeomorphic. In this section we prove the following, which is a restatement of Theorem~\ref{sec:intro:thm:original}, using log transform, knot surgery and the relative genus function. 

\begin{theorem}\label{sec:Producing exotic smooth structures:thm:main} Let $X$ be an arbitrary connected oriented smooth $4$-manifold. Here $X$ is allowed to be closed, or to have $($possibly disconnected$)$ boundary, or to be non-compact. Suppose that $X$ smoothly contains a nucleus. Assume further that $X$ can be embedded into a connected oriented closed smooth $4$-manifold with $b^+_2>1$ whose Seiberg-Witten invariant does not vanish. Then $X$ admits infinitely many exotic smooth structures. 
\end{theorem}
Throughout this section, we fix the notation as follows. 

\begin{definition}\label{sec:construction:def:log_nucleus}
Let $(N,T)$ and $d_T$ be a nucleus and the divisor of $T$, respectively. Let $X$ be an arbitrary connected oriented smooth $4$-manifold which smoothly contains the nucleus $N$. 
Let $Z$ be a closed connected oriented smooth $4$-manifold with $b^+_2>1$ whose Seiberg-Witten invariant does not vanish, and let $\psi:X\to Z$ be a smooth embedding. Assume that $H_2(X;\mathbb{Z})/\textnormal{Tor}$ is a free $\mathbb{Z}$-module, unless otherwise stated. 
\end{definition}
To construct the desired manifolds, we need the following data. 
\begin{definition}[Data set of $X$]\label{sec:cont:def:dataofX}
$(1)$ Let $\widehat{T}, S$ be a basis of $H_2(N;\mathbb{Z})$ such that $[T]=d_T\cdot \widehat{T}$ and $S\cdot \widehat{T}=1$ (see Lemma~\ref{nuclei:lem:basic:log}).
\smallskip
\\
$(2)$ Put $X^0=X-\text{int}\, N$. Since $\partial N$ is a homology $3$-sphere, we get the decomposition 
\begin{equation*}
H_2(X;\mathbb{Z})/\textnormal{Tor}=H_2(N;\mathbb{Z})\oplus (H_2(X^0;\mathbb{Z})/\textnormal{Tor}). 
\end{equation*}
Let 
\begin{equation*}
(\psi\lvert_{X^0})_*:H_2(X^0;\mathbb{Z})/\textnormal{Tor}\to H_2(Z;\mathbb{Z})/\textnormal{Tor}
\end{equation*}
be the homomorphism induced by the embedding $\psi\lvert_{X^0}:X^0\to Z$. 
Since 
\begin{equation*}
(H_2(X^0;\mathbb{Z})/\textnormal{Tor})\big/ \ker (\psi\lvert_{X^0})_*\cong \textnormal{Im}\, (\psi\lvert_{X^0})_*
\end{equation*}
is a free $\mathbb{Z}$-module, we have a direct sum decomposition
\begin{equation*}
H_2(X^0;\mathbb{Z})/\textnormal{Tor}=F_{X^0}\oplus\ker (\psi\lvert_{X^0})_*,
\end{equation*}
where $F_{X^0}$ is a free $\mathbb{Z}$-module isomorphic to $\textnormal{Im}\, (\psi\lvert_{X^0})_*$. Note that $F_{X^0}$ is finitely generated, because $H_2(Z;\mathbb{Z})/\textnormal{Tor}$ is finitely generated. Put $k=\textnormal{rank}\, (F_{X^0})$. Beware that $k$ can be zero. 
\smallskip
\\
$(3)$ 
Let $\{u_1,u_2,\dots,u_k\}$ and $\{ u_{\lambda}\mid \lambda\in \Lambda_1 \}$ be bases of $F_{X^0}$ and $\ker (\psi\lvert_{X^0})_*$, respectively. Note that 
\begin{equation*}
\mathbf{u}:=\{S,\, \widehat{T}, u_1, u_2,\dots,u_k\}\cup \{ u_\lambda \mid \lambda\in \Lambda_1\}
\end{equation*}
is a basis of $H_2(X;\mathbb{Z})/\textnormal{Tor}$. 
\smallskip
\\
$(4)$ For $1\leq i\leq k$, define $g(u_i)$ by 
\begin{equation*}
g(u_i)=\left\{
\begin{array}{ll}
g_{X^0}(u_i), &\text{if $u_i\cdot u_i\geq 0$.} \\
\max \{g_{X^0}(u_i),\, 1\}, &\text{if $u_i\cdot u_i<0$.}
\end{array}
\right.
\end{equation*}
\end{definition} 
\begin{remark}\label{rem:exotic:independent}
When $H_2(X;\mathbb{Z})$ is finitely generated (e.g.\ $X$ is compact), it suffices to set $F_{X^0}=H_2(X^0;\mathbb{Z})/\textnormal{Tor}$, ignoring the above definition. In this case, we do not need to use $\ker (\psi\lvert_{X^0})_*$ in the following argument. 
\end{remark}
\subsection{Construction 1: log transform}\label{subsec:construction:log}
In this subsection, we construct infinitely many distinct exotic smooth structures on $X$ using log transform.
\begin{definition}[Data set of $X_{(p)}$]\label{sec:construction:def:N_i(p)}
$(1)$ Let $p$ be a positive integer which is relatively prime to $d_T$. Let $N_{(p)}$ $(\text{resp}.\ X_{(p)})$ be the $p$-log transform of $N$ $(\text{resp}.\ X)$ along $T$. Let $T_p$ denote the multiple fiber of $N_{(p)}$ (hence, of $X_{(p)}$). When the intersection form of $N$ is even, and $d_T$ is odd, we further assume that $p$ is odd. Note that the intersection forms of $N_{(p)}$ and $N$ are isomorphic by Lemma~\ref{sec:nucleu:lem:property}. 
Lemmas~\ref{nucleus:lem:divisor} and~\ref{sec:nucleu:lem:property} thus imply the existence of a basis $\widehat{T}_p, S_p$ of $H_2(N_{(p)};\mathbb{Z})$ which satisfies 
\begin{equation*}
[T_p]=d_T\cdot \widehat{T}_p,\quad S_p\cdot S_p=S\cdot S,\quad \widehat{T}_{p}\cdot S_{p}=1. 
\end{equation*}
Note $S_1=S$, since $N_{(1)}\cong N$ and $X_{(1)}\cong X$.
\smallskip 
\\
$(2)$ Regard $S_{p}$ as an element of $H_2(X_{(p)};\mathbb{Z})$ through the homomorphism induced by the inclusion. For $p\geq 1$, define $g(S_{p})$ by 
\begin{equation*}
g(S_{p})=\left\{
\begin{array}{ll}
 g_{X_{(p)}}(S_p), &\text{if $S_{p}\cdot S_p\geq 0$.}\\
 \max\{g_{X_{(p)}}(S_p),\, 1\}, &\text{if $S_{p}\cdot S_p<0$.}
\end{array}
\right.
\end{equation*}
\end{definition}

We here define integers $p_n$ $(n=1,2,\dots)$. Roughly speaking, the following conditions require that $p_n$ is sufficiently larger than $p_{n-1}$ for each $n$. 
\begin{definition}\label{sec:construction::def:p_n}Put $p_1=1$. Define an integer sequence $p_n$ $(n\geq 2)$ so that it satisfies the following conditions (i)--(iv).  
\begin{itemize}
 \item [(i)] $p_{n}>p_{n-1}$, \, for each $n\geq 2$.
 \item [(ii)] $d_T(p_2-1)+u_i\cdot u_i>2g(u_i)-2$, \, for each $1\leq i\leq k$.
\item [(iii)] $d_T(p_n-1)+S\cdot S>2g(S_{p_{n-1}})-2$, \, for each $n\geq 2$.
 \item[(iv)] Each $p_n$ $(n\geq 2)$ is relatively prime to $d_T$. 
\end{itemize}
In the case where the intersection form of $N$ is even (i.e. $S\cdot S$ is even), and $d_T$ is odd, we further assume the following (v). 
\begin{itemize}
\item [(v)] Every $p_n$ $(n\geq 2)$ is odd. 
\end{itemize}
\end{definition}

We are now ready to define the desired manifolds. 
\begin{definition}\label{sec:infinite_compact:def:X_n(n>0)}
Let $X_n$ (resp.\ $Z_n$) $(n\geq 1)$ be the $p_n$-log transform of $X$ (resp.\ $Z$) along the torus $T$ in the nucleus $N$. 
\end{definition}

Note that $X_1$ (resp.\ $Z_1$) is diffeomorphic to $X$ (resp.\ $Z$). We here check topological types of $X_n$'s. 

\begin{lemma}\label{sec:construction:lem:homeo}Every $X_n$ $(n\geq 1)$ is homeomorphic to $X$. 
\end{lemma}
\begin{proof}Since $X_n$ is obtained from $X$ by replacing the copy of $N$ with $N_{(p_n)}$, Lemma \ref{sec:nucleu:lem:property} shows the claim. 
\end{proof}
\subsection{Detecting smooth structures}\label{subsec:detect}
In this subsection, we detect smooth structures of $X_n$'s obtained in Subsection~\ref{subsec:construction:log}. To see differences of the relative genus functions, we choose an element $(Q, \mathbf{d}, \mathbf{g})$ as follows.
\begin{definition}\label{def:Q,d,g}Let $\Lambda$ be the disjoint union $\{0,1,\dots,k+1\}\coprod \Lambda_1$, and put $\lambda_0=0\in \Lambda$. 
Let $Q\in \textnormal{Sym}_\Lambda(\mathbb{Z})$ be the intersection matrix of $H_2(X;\mathbb{Z})/ \textnormal{Tor}$ given by the basis $\mathbf{u}$. Define $\mathbf{d}=(d_\lambda)_{\lambda\in \Lambda}$ and $\mathbf{g}=(g_\lambda)_{\lambda\in \Lambda-\{\lambda_0\}}$ as follows:
\begin{align*}
d_0=1, \quad d_1&=d_T, \quad d_i=1 \; (2\leq i\leq k+1), \quad d_\lambda=0 \; (\lambda\in \Lambda_1), \\
g_1&=1,\quad g_i=g(u_{i-1}) \; (2\leq i\leq k+1), \quad g_\lambda=0 \; (\lambda\in \Lambda_1). 
\end{align*}

\end{definition}
For this $(Q, \mathbf{d}, \mathbf{g})$, we obtain the following evaluation.
\begin{proposition}\label{detect:genus:log}Let $(Q, \mathbf{d}, \mathbf{g})\in \textnormal{Sym}_\Lambda(\mathbb{Z})\times \mathbb{Z}^{\Lambda}\times \mathbb{Z}^{\Lambda-\{\lambda_0\}}$ be the one in Definition~\ref{def:Q,d,g}. Then the following inequalities hold. 
\begin{itemize}
 \item For each $n\geq 1$, \, $G_{X_n}(Q, \mathbf{d}, \mathbf{g})\leq g(S_{p_n})$.
 \item For each $n\geq 2$, \, $g(S_{p_{n-1}})< G_{X_n}(Q, \mathbf{d}, \mathbf{g})$.
\end{itemize}
Consequently $G_{X_n}\neq G_{X_{m}}$ for any $n\neq m$. 
\end{proposition}
\begin{proof} Put $\mathbf{u}_{n}=\{S_{p_n}, \widehat{T}_{p_n}, u_1,u_2,\dots,u_k\}\cup \{ u_\lambda \mid \lambda\in \Lambda_1\}$, and regard $\mathbf{u}_{n}$ as the subset of $H_2(X_n;\mathbb{Z})/ \textnormal{Tor}$ through the natural inclusions $N_{(p_n)}\to X_n$ and $X^0\to X_n$. 
This $\mathbf{u}_{n}$ is clearly a basis of $H_2(X_n;\mathbb{Z})/ \textnormal{Tor}$. By Definitions~\ref{sec:cont:def:dataofX} and \ref{sec:construction:def:N_i(p)}, we have 
\begin{equation*}
\text{$G_{X_n}(Q, \mathbf{d}, \mathbf{g})\leq G_{X_n, Q, \mathbf{d}, \mathbf{g}}(\mathbf{u}_{n})\leq g(S_{p_n})$, for each $n\geq 1$.} 
\end{equation*}

We next show $G_{X_n}(Q, \mathbf{d}, \mathbf{g}) > g(S_{p_{n-1}})$ for each $n\geq 2$. Here note that every $u_\lambda$ $(\lambda\in \Lambda_1)$ satisfies $u_\lambda \cdot w=0$ for any $w\in H_2(X;\mathbb{Z})/ \textnormal{Tor}$, because $u_\lambda\in \ker (\psi\lvert_{X^0})_*$. Suppose $G_{X_n}(Q, \mathbf{d}, \mathbf{g}) \leq g(S_{p_{n-1}})$. Then there exists a basis $\{v_\lambda\mid\lambda\in\Lambda\}$ of $H_2(X_n;\mathbb{Z})/ \textnormal{Tor}$ satisfying the following (i)--(v). 
\begin{enumerate}
 \item [(i)] $v_0\cdot v_0=S\cdot S$,\quad $v_1\cdot v_1=0$,\quad $v_{i}\cdot v_{i}=u_{i-1}\cdot u_{i-1}$ $(2\leq i\leq k+1)$.
 \item [(ii)] $v_0$ is represented by a smoothly embedded closed oriented surface of genus $g(S_{p_{n-1}})$.
 \item [(iii)] $d_Tv_1$ is represented by a smoothly embedded torus. 
 \item [(iv)] Each $v_i$ $(2\leq i\leq k+1)$ is represented by a smoothly embedded closed oriented surface of genus $g(u_{i-1})$. 
 \item [(v)] For any $\lambda\in \Lambda_1$ and any $w\in H_2(X_n;\mathbb{Z})/ \textnormal{Tor}$, they satisfy $v_\lambda \cdot w=0$.
\end{enumerate}
Since $\mathbf{u}_{n}$ is a basis, for each $\lambda\in \Lambda$, there exists a set $\{a^{(\lambda)}_\mu \in \mathbb{Z}\mid \mu\in \Lambda\}$ satisfying
\begin{equation*}
v_\lambda=a^{(\lambda)}_0S_{p_n}+a^{(\lambda)}_1\widehat{T}_{p_n}+a^{(\lambda)}_2u_1+a^{(\lambda)}_3u_2+\dots+a^{(\lambda)}_{k+1}u_{k}+\sum_{\mu\in \Lambda_1}a^{(\lambda)}_\mu u_{\mu}, 
\end{equation*}
where $a^{(\lambda)}_\mu=0$ all but finitely many $\mu$'s. 
We may regard each $v_\lambda$ $(\lambda\in \Lambda)$ as an element of $H_2(Z_{n};\mathbb{Z})$ (up to torsion) through the inclusion $X_n\to Z_{n}$. 

Theorem~\ref{Thm:SW:log} gives a class $K$ of $H^2(Z_n;\mathbb{Z})$ such that
\begin{equation*}
L^+:=K+(p_n-1)PD([T_{p_n}])\quad \text{and}\quad L^-:=K-(p_n-1)PD([T_{p_n}])
\end{equation*}
are Seiberg-Witten basic classes of $Z_n$. 
Applying the adjunction inequality to $v_{0}$ and the basic classes $L^+$ and $L^-$, we get the inequalities below.  
\begin{align*}
&2g(S_{p_{n-1}})-2 \geq \left|\langle K, v_0\rangle + d_T(p_n-1)a^{(0)}_0 \right|+S\cdot S, \\
&2g(S_{p_{n-1}})-2 \geq \left|\langle K, v_0\rangle -d_T(p_n-1)a^{(0)}_0 \right|+S\cdot S. 
\end{align*}
These inequalities imply
\begin{equation*}
2g(S_{p_{n-1}})-2\geq d_T(p_n-1)\left| a^{(0)}_{0}\right|+S\cdot S. 
\end{equation*}
The condition (iii) of $p_n$ in Definition~\ref{sec:construction::def:p_n} thus shows $a^{(0)}_{0}=0$. We next apply the adjunction inequality to $v_1$ and the basic classes $L^+$ and $L^-$. We then get 
 \begin{align*}
0&=\left|\langle K, v_1\rangle + d_T(p_n-1)a^{(1)}_0 \right|=\left|\langle K, v_1\rangle  -d_T(p_n-1)a^{(1)}_0  \right|. 
\end{align*}
Since $d_T$ and $p_n-1$ are positive, this equality shows $a^{(1)}_{0}=0$. We finally apply the adjunction inequality to each $v_{i}$ $(2\leq i\leq k+1)$ and the basic classes $L^+$ and $L^-$. We then get the inequalities below. 
\begin{align*}
&2g(u_{i-1})-2 \geq \left|\langle K, v_i\rangle + d_T(p_n-1)a^{(i)}_0 \right|+u_{i-1}\cdot u_{i-1}, \\
&2g(u_{i-1})-2 \geq \left|\langle K, v_i\rangle -d_T(p_n-1)a^{(i)}_0 \right|+u_{i-1}\cdot u_{i-1}. 
\end{align*}
These inequalities give
\begin{equation*}
2g(u_{i-1})-2\geq d_T(p_n-1)\left|a^{(i)}_0\right|+u_{i-1}\cdot u_{i-1}. 
\end{equation*}
The conditions (i) and (ii) of $p_n$ in Definition~\ref{sec:construction::def:p_n} thus show $a^{(i)}_{0}=0$ $(2\leq i\leq k+1)$. 

The condition (v) of the basis $\{v_\lambda \mid \lambda \in \Lambda \}$ shows $a^{(\lambda)}_0=0$ for any $\lambda\in \Lambda_1$. 
We now have $a^{(\lambda)}_0=0$ for all $\lambda\in \Lambda$. This contradicts the assumption that $\{v_\lambda \mid \lambda \in \Lambda \}$ is a basis of $H_2(X_n;\mathbb{Z})/\textnormal{Tor}$. 
\end{proof}


This proposition and Lemma~\ref{sec:construction:lem:homeo} clearly give the following. 
\begin{theorem}\label{sec:detect:thm:diffeo_type}
Assume that $H_2(X;\mathbb{Z})/\textnormal{Tor}$ is a free $\mathbb{Z}$-module. Then the smooth $4$-manifolds $X_n$ $(n\geq 1)$ are all homeomorphic to $X$ but mutually non-diffeomorphic with respect to the given orientations. 
\end{theorem}
We are now ready to prove Theorem~\ref{sec:Producing exotic smooth structures:thm:main}. 
\begin{proof}[Proof of Theorem~\ref{sec:Producing exotic smooth structures:thm:main}] When $H_2(X;\mathbb{Z})/\textnormal{Tor}$ is a free module, Theorem~\ref{sec:Producing exotic smooth structures:thm:main} clearly follows from the theorem above. We prove the general case, that is, we allow $H_2(X;\mathbb{Z})/\textnormal{Tor}$ to be non-free module. We use the same definitions as in Definitions~\ref{sec:construction:def:log_nucleus}, \ref{sec:cont:def:dataofX}, \ref{sec:construction:def:N_i(p)}, \ref{sec:construction::def:p_n} and \ref{sec:infinite_compact:def:X_n(n>0)}, except for the following $(\psi\lvert_{X^0})_*$ and $F_{X^0}$. Let 
\begin{equation*}
(\psi\lvert_{X^0})_*:H_2(X^0;\mathbb{Z})\to H_2(Z;\mathbb{Z})/\textnormal{Tor}
\end{equation*}
 be the homomorphism induced by the embedding $\psi\lvert_{X^0}:X^0\to Z$. Since 
\begin{equation*}
H_2(X^0;\mathbb{Z})/ \ker (\psi\lvert_{X^0})_* \cong \textnormal{Im}\, (\psi\lvert_{X^0})_*
\end{equation*}
is a free $\mathbb{Z}$-module, $H_2(X^0;\mathbb{Z})$ has a direct sum decomposition 
\begin{equation*}
H_2(X^0;\mathbb{Z})=F_{X^0}\oplus \ker (\psi\lvert_{X^0})_*, 
\end{equation*}
where $F_{X^0}$ is a finitely generated free $\mathbb{Z}$-module isomorphic to $\textnormal{Im}\, (\psi\lvert_{X^0})_*$. Put
\begin{equation*}
F_0=H_2(N;\mathbb{Z})\oplus F_{X^0}\quad \text{and}\quad F_1=\ker (\psi\lvert_{X^0})_*.
\end{equation*}
Similarly to the proof of Proposition~\ref{detect:genus:log}, we can show $G_{X_n,F_0, F_1}\neq G_{X_m,F_0, F_1}$ for any $n\neq m$ (Consider the subset $\mathbf{u}_{n}'=\{S_{p_n}, \widehat{T}_{p_n}, u_1,u_2,\dots,u_k\}\cup \ker (\psi\lvert_{X^0})_*$ of $H_2(X_n;\mathbb{Z})$. Note that $\mathbf{u}_{n}'$ is a generating set of $H_2(X;\mathbb{Z})$). Therefore the desired claim follows. 
\end{proof}
\subsection{Construction 2: knot surgery}
In this subsection, we construct infinitely many distinct exotic smooth structures on $X$ using knot surgery. 

Recall Definitions~\ref{sec:construction:def:log_nucleus} and \ref{sec:cont:def:dataofX}. In this subsection, as an additional condition of the nucleus $(N,T)$, we assume $d_T=1$. 
\begin{definition}\label{knot surgery:def:g(S_K)}$(1)$ Let $K$ be a knot in $S^3$, and let $\deg (\triangle_K)$ denote the maximal degree of the symmetrized Alexander polynomial of $K$. Let $S_K$ be the primitive class of $H_2(N_K;\mathbb{Z})$ which satisfies $[T]\cdot S_K=1$ and $S_K\cdot S_K=S\cdot S$ (see Lemma~\ref{lem:basic:nuclei:knot} and its proof).
\\
$(2)$ Regard $S_K$ as an element of $H_2(X_K;\mathbb{Z})$ through the homomorphism induced by the inclusion. Define $g(S_{K})$ by
\begin{equation*}
g(S_K)=\left\{
\begin{array}{ll}
g_{X_K}(S_K), &\text{if $S\cdot S\geq 0$.}\\
\max\{g_{X_K}(S_K),\, 1\}, &\text{if $S\cdot S<0$.}
\end{array}
\right.
\end{equation*}
\end{definition}

We here define an infinite sequence $K_n$ $(n=1,2,\dots)$ of knots in $S^3$. Roughly speaking, the following conditions require that $\deg (\triangle_{K_n})$ is sufficiently larger than $\deg (\triangle_{K_{n-1}})$ for each $n$ (Such a sequence of knots is known to exist (cf.\ p.366 of \cite{FS2})). 
\begin{definition}\label{def:K_n}Let $K_1$ be the unknot in $S^3$. 
Let $K_n$ $(n\geq 2)$ be an infinite sequence of knots in $S^3$ which satisfies the following conditions (i)--(iii).  
\begin{itemize}
 \item [(i)] $\deg (\triangle_{K_{n}})>\deg (\triangle_{K_{n-1}})$, \, for each $n\geq 2$.
 \item [(ii)] $2\deg (\triangle_{K_2})+u_i\cdot u_i>2g(u_i)-2$, \, for each $1\leq i\leq k$.
\item [(iii)] $2\deg (\triangle_{K_n})+S\cdot S>2g(S_{K_{n-1}})-2$, \, for each $n\geq 2$.
\end{itemize}
\end{definition}
We are now ready to define the desired manifolds. 
\begin{definition}\label{sec:infinite_compact:def:knot surgery}
Let $Y_n$ (resp.\ $Z_n$) $(n\geq 1)$ be the knot surgery of $X$ (resp.\ $Z$) along the torus $T$ with the knot $K_n$. 
\end{definition}
Note that $Y_1$ (resp.\ $Z_1$) is diffeomorphic to $X$ (resp.\ $Z$). The reason is as follows. Since $K_1$ is unknot, $S^3-\nu(K_1)$ is diffeomorphic to $S^1\times D^2$. Therefore $Y_1$ is the $1$-log transform of $X$ on the cusp neighborhood (cf.\ Remark~\ref{rem:multiple fiber}). Hence $Y_1\cong X_{(1)}\cong X$ (cf.\ \cite{GS}). 

Similarly to the proof of Proposition~\ref{detect:genus:log}, we obtain the following evaluation. 
\begin{proposition}\label{infinite:prop:knot}Let $(Q, \mathbf{d}, \mathbf{g})$ be the one in Definition~\ref{def:Q,d,g} $($Note $d_T=1$ in this case$)$. Then the following inequalities hold.
\begin{itemize}
 \item For each $n\geq 1$, \,$G_{Y_n}(Q, \mathbf{d}, \mathbf{g})\leq g(S_{K_n})$.
 \item For each $n\geq 2$, \, $g(S_{K_{n-1}})< G_{Y_n}(Q, \mathbf{d}, \mathbf{g})$.
\end{itemize}
Consequently $G_{Y_n}\neq G_{Y_{m}}$ for any $n\neq m$. 
\end{proposition}
This proposition and Lemma~\ref{lem:basic:nuclei:knot} clearly provide us the following. 
\begin{theorem}\label{sec:construct:knot:detect}Assume that $H_2(X;\mathbb{Z})/\textnormal{Tor}$ is a free $\mathbb{Z}$-module and that $d_T=1$. Then the smooth $4$-manifolds $Y_n$ $(n\geq 1)$ are all homeomorphic to $X$, but mutually non-diffeomorphic with respect to the given orientations. 
\end{theorem}
\begin{remark}Similarly to the proof of Theorem~\ref{sec:Producing exotic smooth structures:thm:main}, we can prove Theorem~\ref{sec:construct:knot:detect} even when $H_2(X;\mathbb{Z})/\textnormal{Tor}$ is not a free $\mathbb{Z}$-module. Note that, in this case, we need to slightly modify the construction of $Y_n$'s, similarly to the proof of Theorem~\ref{sec:Producing exotic smooth structures:thm:main}. 
\end{remark}
\subsection{Strengthening the constructions}
In the previous subsections, we did not exclude the possibility that some of the $4$-manifolds are orientation reversing diffeomorphic. In this subsection, we exclude this possibility by restricting the conditions of $p_n$'s and $K_n$'s in Definitions~\ref{sec:construction::def:p_n} and \ref{def:K_n}. 

Let $(N,T)$, $X$, $Z$ be as in Definition~\ref{sec:construction:def:log_nucleus}. 
Recall the definitions of $X_n$'s and $Y_n$'s in Definitions~\ref{sec:infinite_compact:def:X_n(n>0)} and \ref{sec:infinite_compact:def:knot surgery}. We here strengthen the constructions of these manifolds as follows. 
\begin{definition}\label{subsec:strength:def:p_n}
Assume that the integer sequence $p_n$ $(n\geq 1)$ in Definition~\ref{sec:construction::def:p_n} further satisfies the following (vi) and (vii). Then define $\widehat{X}_n=X_n$ $(n\geq 1)$. 
\begin{itemize}
\item [(vi)] $d_T(p_2-1)-u_i\cdot u_i>2g(u_i)-2$, \, for each $1\leq i\leq k$.
\item [(vii)] $d_T(p_n-1)-S\cdot S>2g(S_{p_{n-1}})-2$, \, for each $n\geq 2$.
\end{itemize}
\end{definition}

\begin{definition}\label{subsec:def:knot:strength}Assume that $d_T=1$ and that the sequence $K_n$ $(n\geq 2)$ of knots in Definition~\ref{def:K_n} further 
satisfies the following conditions (iv) and (v). Then define $\widehat{Y}_n=Y_n$ $(n\geq 1)$.
\begin{itemize}
\item [(iv)] $2\deg (\triangle_{K_2})-u_i\cdot u_i>2g(u_i)-2$, \, for each $1\leq i\leq k$. 
\item [(v)] $2\deg (\triangle_{K_n})-S\cdot S>2g(S_{K_{n-1}})-2$, \, for each $n\geq 2$. 
\end{itemize}
\end{definition}

 Similarly to the proof of Proposition~\ref{detect:genus:log}, we obtain the following stronger evaluations.
\begin{proposition}\label{detect:genus:log:strengthen}Let $(Q, \mathbf{d}, \mathbf{g})$ be the one in Definition~\ref{def:Q,d,g}. Then the following inequalities hold. 
\begin{itemize}
 \item For each $n\geq 1$, \, $G^{\pm}_{\widehat{X}_n}(Q, \mathbf{d}, \mathbf{g})\leq g(S_{p_n})$.
 \item For each $n\geq 2$, \, $g(S_{p_{n-1}})< G^{\pm}_{\widehat{X}_n}(Q, \mathbf{d}, \mathbf{g})$.
 \end{itemize}
Consequently $G^{\pm}_{\widehat{X}_n}\neq G^{\pm}_{\widehat{X}_{m}}$ for any $n\neq m$. 
\end{proposition}
\begin{proposition}\label{infinite:prop:knot:strengthen}Let $(Q, \mathbf{d}, \mathbf{g})$ be the one in Definition~\ref{def:Q,d,g} $($Note $d_T=1$ in this case$)$. Then the following inequalities hold. 
\begin{itemize}
 \item For each $n\geq 1$, \, $G^{\pm}_{\widehat{Y}_n}(Q, \mathbf{d}, \mathbf{g})\leq g(S_{K_n})$.
 \item For each $n\geq 2$, \, $g(S_{K_{n-1}})< G^{\pm}_{\widehat{Y}_n}(Q, \mathbf{d}, \mathbf{g})$.
 \end{itemize}
Consequently $G^{\pm}_{\widehat{Y}_n}\neq G^{\pm}_{\widehat{Y}_{m}}$ for any $n\neq m$. 
\end{proposition}

These propositions provide us the stronger versions of Theorems~\ref{sec:detect:thm:diffeo_type} and \ref{sec:construct:knot:detect}. 

\begin{theorem}\label{sec:construct:knot:detect:strengthen}Assume that $H_2(X;\mathbb{Z})/\textnormal{Tor}$ is a free $\mathbb{Z}$-module. 
Then the following hold.\\
$(1)$ $\widehat{X}_n$ $(n\geq 1)$ are all homeomorphic to $X$, but mutually non-diffeomorphic for any orientations. \\
$(2)$ $\widehat{Y}_n$ $(n\geq 1)$ are all homeomorphic to $X$, but mutually non-diffeomorphic for any orientations. 
\end{theorem}

\subsection{Remarks}\label{subsection:remarks}We conclude this section making some remarks. 
\smallskip

(1) Modifying the constructions of $\widehat{X}_n$'s and $\widehat{Y}_n$'s similarly to the proof of Theorem~\ref{sec:Producing exotic smooth structures:thm:main}, we can show that Theorem~\ref{sec:construct:knot:detect:strengthen} holds even when $H_2(X;\mathbb{Z})/\textnormal{Tor}$ is not a free $\mathbb{Z}$-module.
\smallskip

(2) The constructions of $X_n$'s, $Y_n$'s, $\widehat{X}_n$'s and $\widehat{Y}_n$'s depend on the choice of the embedding $\psi: X\to Z$. When $H_2(X;\mathbb{Z})$ is finitely generated (e.g.\ $X$ is compact)Cwe can arrange the constructions so that they do not depend on the choice of $\psi$. See Remark~\ref{rem:exotic:independent}.
\smallskip


(3) While we required the closed 4-manifold $Z\supset X$ to be $b_2^+(Z)>1$ in Theorem~\ref{sec:Producing exotic smooth structures:thm:main}, the same claim also holds if $\pi_1(Z)\cong 1$, $b_2^+(Z)=1$ and $b_2^-(Z)\leq 9$. For the log transform construction, this can be seen as follows (the knot surgery case is similar). 
The adjunction inequality for the above $Z$ was given, for example, in \cite{P1.5}, though the inequality was proved only for surfaces with the non-negative self-intersection numbers. We can thus prove the desired claim by modifying the conditions of the integer sequence $p_n$ $(n\geq 1)$ and the proof of Proposition~\ref{detect:genus:log}, namely as follows. Choose the class $S$ so that $S\cdot S>0$. Then consider a basis $\mathbf{u}_{n}'=\{S_{p_n}, \widehat{T}_{p_n}, u_1+l_1S_{p_n}, u_2+l_2S_{p_n},\dots,u_k+l_{k}S_{p_n}\}\cup \{ u_\lambda\mid \lambda\in \Lambda_1\}$ of $H_2(X_n;\mathbb{Z})/ \textnormal{Tor}$, where each $l_i$ is an integer satisfying $(u_i+l_iS_{p_n})^2\geq 0$. Let $Q\in \textnormal{Sym}_\Lambda(\mathbb{Z})$ be the intersection matrix given by the basis $\mathbf{u}_{n}'$, and let $\mathbf{d}\in \mathbb{Z}^\Lambda$ be the one in Definition~\ref{def:Q,d,g}. Since $g_{X_n}(u_{i}+l_{i}S_{p_n})\leq g_{X^0}(u_{i})+g_{N_{(p_n)}}(l_{i}S_{p_n})$, the rest of the argument is similar, except that we use infinitely many different $\mathbf{g}\in \mathbb{Z}^{\Lambda-\{\lambda_0\}}$ to see differences of relative genus functions of $X_n$'s $(n\geq 1)$. 

\section{Compact Stein $4$-manifolds}\label{sec:stein}In this section, we briefly recall Stein $4$-manifolds. For the definition of basic terms and more details, the reader can consult \cite{GS} and \cite{OS1}. 
In this paper, we use Seifert framings and sometimes abbreviate them to framings. (When a knot goes over $4$-dimensional $1$-handles, then convert the diagram into the dotted circle notation and calculate its Seifert framing (cf.~\cite{GS}).) We use the following terminologies in the rest of this paper. 

\begin{definition}\label{sec:stein:def:handlebody}
$(1)$ For a Legendrian knot $K$ in $\#n(S^1\times S^2)$ $(n\geq 0)$, we denote by $tb(K)$ and $r(K)$ the Thurston-Bennequin number and the rotation number of $K$, respectively.
\\
$(2)$ We call a compact connected oriented $4$-dimensional handlebody a \textit{$2$-handlebody} if it consists of one $0$-handle and $1$- and $2$-handles. 
\\
$(3)$ We call a $2$-handlebody a \textit{Legendrian handlebody} if its $2$-handles are attached along an oriented framed Legendrian link in $\partial(D^4\cup \textit{$1$-handles})=\#n(S^1\times S^2)$ $(n\geq 0)$. It is known that every $2$-handlebody can be changed into a Legendrian handlebody by an isotopy of the attaching link of $2$-handles, and orienting its components.
\\
$(4)$ We call a Legendrian handlebody a \textit{Stein handlebody} if the framing of its each $2$-handle $K$ is $tb(K)-1$. 
\end{definition}
\begin{remark}
For a given Legendrian knot, we can decrease its Thurston-Bennequin number by an arbitrary positive integer, by locally adding ``zig-zags'' to the knot. 
\end{remark}

Next we recall the following useful theorem. 

\begin{theorem}[Eliashberg~\cite{E1}, Gompf~\cite{G}. cf.\ \cite{GS}]\label{sec:stein:th:stein_existence}
A compact, connected, oriented, smooth $4$-manifold admits a Stein structure if and only if it can be represented as a Stein handlebody. 
\end{theorem}
We call a compact smooth $4$-manifold with a Stein structure a \textit{compact Stein $4$-manifold}. 

\begin{figure}[ht!]
\begin{center}
\includegraphics[width=1.3in]{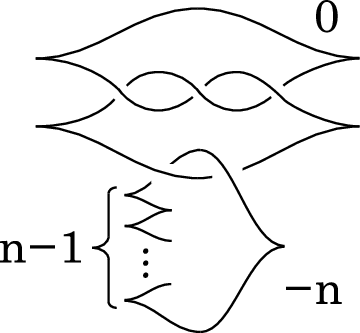}
\caption{A Stein handlebody decomposition of $G(n)$ $(n\geq 2)$.}
\label{fig:nuclei_fig4}
\end{center}
\end{figure}

\begin{example}\label{ex:gompf:stein}Gompf nucleus $G(n)$ $(n\geq 2)$ admits a Stein handlebody decomposition as shown in Figure~\ref{fig:nuclei_fig4}. 
\end{example}
One can similarly construct many other nuclei which admit Stein structures. These Stein nuclei are useful to construct exotic smooth structures. See Section~\ref{sec:stein and exotic}. 

Compact Stein $4$-manifolds are known to admit the following useful embeddings, where minimal means that there exists no smoothly embedded $2$-sphere with the self-intersection number $-1$. 
\begin{theorem}[Lisca-Mati\'{c}~\cite{LM1}]\label{sec:stein:th:closing of Stein:complex}
Every compact Stein $4$-manifold can be embedded into a minimal closed complex surface of general type with $b_2^+>1$. 
\end{theorem}
\begin{theorem}[Akbulut-Ozbagci~\cite{AO3}]\label{sec:stein:th:closing of Stein:symplectic}
Every compact Stein $4$-manifold can be embedded into a simply connected, minimal, closed, symplectic $4$-manifold with $b_2^+>1$.
\end{theorem}
\begin{proof}(simply connectedness). Since ``simply connected'' is not claimed in \cite{AO3}, we explain only this part for completeness. We follow the proof in~\cite{AO3}. We first attach $2$-handles to a given Stein $4$-manifold to make it simply connected Stein $4$-manifold.  Then we apply the procedure prescribed in~\cite{AO3}. Namely, we attach $2$-handles to this simply connected Stein $4$-manifold so that it admits a Lefschetz fibration over $D^2$. As their construction shows, this fibration can be extended to a Lefschetz fibration $X\to S^2$ by attaching $F\times D^2$ (i.e.\ $2$-, $3$- and $4$-handles), where $F$ denotes a regular fiber. As shown in~\cite{AO3}, $X$ is a (possibly non-minimal) closed symplectic $4$-manifold with $b_2^+>1$. Note that $X$ and $X-\nu (F)$ is simply connected. Taking the fiber sum of two copies of $X$, we get a simply connected, minimal, closed, symplectic $4$-manifold with $b_2^+>1$ by a theorem of Stipsicz~\cite{Sti}. 
\end{proof}
\begin{remark}\label{rem:SW of general type}Every closed symplectic $4$-manifold with $b_2^+>1$ is known to have the non-vanishing Seiberg-Witten invariant (cf.\ \cite{GS}). Furthermore, it is known that every minimal closed complex surface of general type with $b_2^+>1$ has only one Seiberg-Witten basic class up to sign and that the square of the basic class is positive (see Theorems 3.4.22 and 3.4.19 of \cite{GS}). 
\end{remark}

\section{Corks and $W$-modifications}\label{sec:cork}
In this section, we briefly recall corks~(\cite{A1}, \cite{AY1}) and $W$-modifications (\cite{AY5}). 
\subsection{Corks}
\begin{definition}
Let $C$ be a compact contractible (Stein) $4$-manifold with boundary and $\tau: \partial C\to \partial C$ an involution on the boundary. 
We call $(C, \tau)$ a \textit{cork} if $\tau$ extends to a self-homeomorphism of $C$, but cannot extend to any self-diffeomorphism of $C$. 
For a cork $(C,\tau)$ and a smooth $4$-manifold $X$ which contains $C$, a {\it cork twist of $X$ along $(C,\tau)$ }is defined to be the smooth $4$-manifold obtained from $X$ by removing the submanifold $C$ and regluing it via the involution $\tau$.   
\end{definition}
Note that Boyer's theorem~\cite{B} tells that every self-diffeomorphism of the boundary $\partial C$ extends to a self-homeomorphism of $C$, in the case where $C$ is a compact contractible smooth $4$-manifold. 

\begin{definition}Let $W_n$ be the compact smooth $4$-manifold in Figure~$\ref{fig:nuclei_fig5}$. Let $f_n:\partial W_n\to \partial W_n$ be the diffeomorphism obtained by first surgering $S^1\times D^3$ to $D^2\times S^2$ in the interior of $W_n$, then surgering the other embedded $D^2\times S^2$ back to $S^1\times D^3$ (i.e. replacing the dot and 0 in Figure ~$\ref{fig:nuclei_fig5}$). 
\begin{figure}[ht!]
\begin{center}
\includegraphics[width=1.1in]{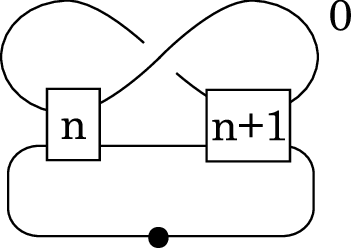}
\caption{$W_n$}
\label{fig:nuclei_fig5}
\end{center}
\end{figure}
\end{definition}We can check that $W_n$ is contractible and that $f_n$ is an involution as follows. 
Changing the crossing of the attaching circle of the $2$-handle of $W_n$ as in the left picture of Figure~\ref{fig:nuclei_fig6}, we get $W_{n-1}$ by isotopy (\cite{AY1}). Repeating this process (i.e.\ by homotopy), we get $W_0(\cong D^4)$. Hence $W_n$ is contractible. Since the diagram of $W_n$ is induced from a symmetric link, the $4$-manifold obtained by replacing the dot and 0 is diffeomorphic to $W_n$. Fixing the identification of $W_n$ (hence, of $\partial W_n$) by Figure~$\ref{fig:nuclei_fig5}$, the diffeomorphism $f_n$ becomes an involution of $\partial W_n$.

\begin{figure}[ht!]
\begin{center}
\includegraphics[width=2.9in]{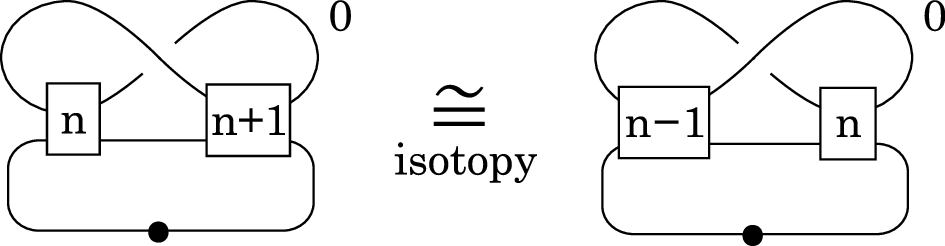}
\caption{$W_{n-1}$}
\label{fig:nuclei_fig6}
\end{center}
\end{figure}


\begin{theorem}[Theorem 2.5 of \cite{AY1}]\label{th:cork}
For $n\geq 1$, the pair $(W_n, f_n)$ is a cork. 
\end{theorem}
\subsection{$W$-modifications}
In this subsection, we recall $W$-modifications introduced in \cite{AY5}. We first define them for smooth $2$-handlebodies and later redefine them for Legendrian handlebodies. In this paper, the words the ``attaching circle of a $2$-handle'' and a ``smoothly embedded surface'' are often abbreviated to  a ``$2$-handle'' and a ``surface'', if they are clear from the context. 
\begin{definition}Let $p$ be a positive integer, and let $K$ be a $2$-handle of a (smooth) $2$-handlebody $X$. Take a small segment of the attaching circle of $K$ as in the first row of Figure~\ref{fig:nuclei_fig7}.

We call the local operations shown in the left and the right side of Figure~\ref{fig:nuclei_fig7} a \textit{$W^+(p)$-modification} to $K$ and a \textit{$W^-(p)$-modification} to $K$, respectively. Here we do not change the framing of $K$ (ignore the orientations shown in the figure). They are clearly related by a cork twist along $(W_1,f_1)$ as shown in the figure.

 We will call the $0$-framed $2$-handle $\gamma$ on the left (or  right) side of the Figure~\ref{fig:nuclei_fig7} the \textit{auxiliary $2$-handle} of the $W^{\pm}(p)$-modification of $K$. We will use the same symbol $K$ for the new $2$-handle obtained from the original $K$ of $X$ by the modification.

For convenience, we refer the $W^+(0)$- and $W^-(0)$-modifications as undone operations. 
For brevity, sometimes we will call these operations  $W^+$- and $W^-$-modifications when we do not need to specify the coefficients, or call them as $W$-modifications when we do not need to specify both the coefficient and $\pm$. 
\begin{figure}[ht!]
\begin{center}
\includegraphics[width=4.5in]{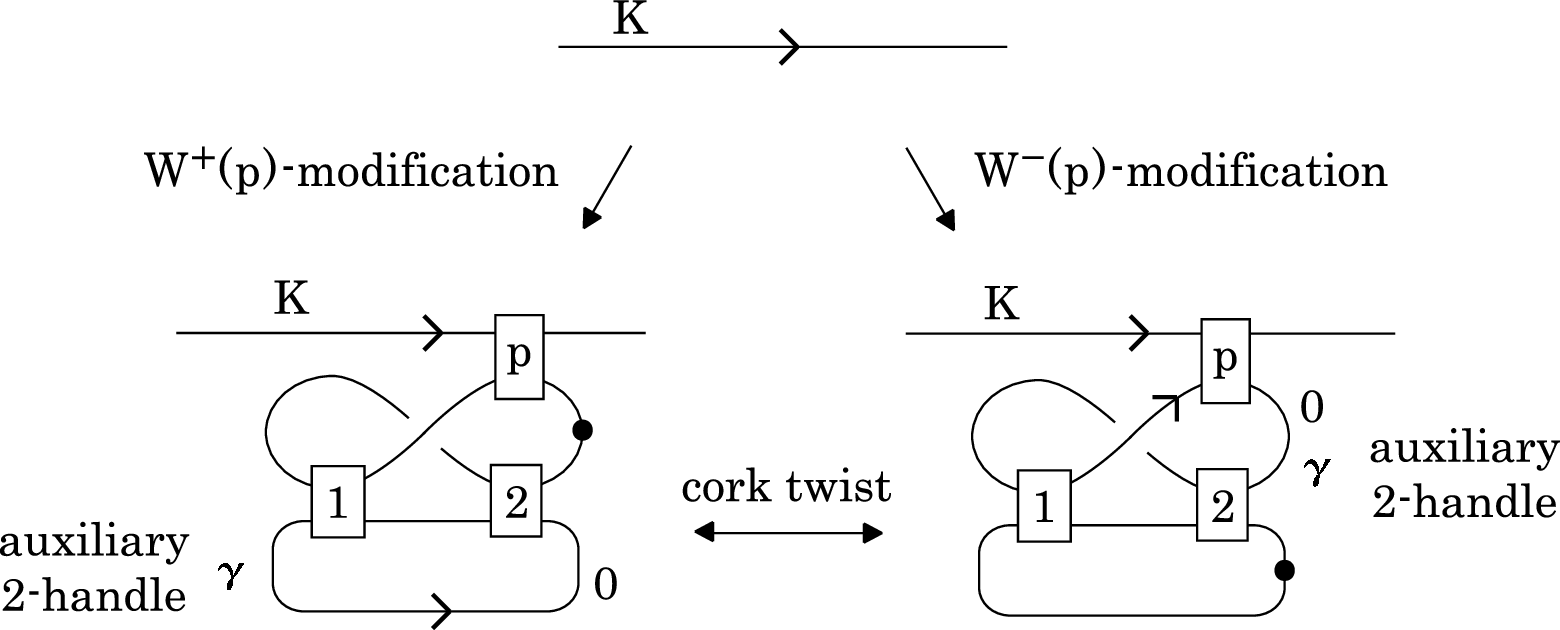}
\caption{$W^{\pm}(p)$-modifications $(p\geq 1)$  (the framing of $K$ is unchanged)}
\label{fig:nuclei_fig7}
\end{center}
\end{figure}
\end{definition}
Similarly we can talk about  \textit{$W^{\pm}(p)$-modification} for any cork $(W,f)$ coming from a symmetric link. However, for simplicity,  we will discuss the effects of $W$-modifications only when $(W,f)=(W_1,f_1)$. 

It is easy to check the following lemmas. 

\begin{lemma}[Proposition~4.2 of \cite{AY5}]\label{prop:relation of attachments} Let $X$ be a $2$-handlebody. Then any $W$-modification to $X$ does not change the isomorphism classes of the fundamental group, the integral homology groups, the integral homology groups of the boundary $\partial X$, and the intersection form of $X$. 
\end{lemma}
\begin{lemma}[Proposition~4.5 of \cite{AY5}]\label{lem:W:embedding}
Let $X$ be a $2$-handlebody. Let $X^+$ $($resp.\ $X^-$$)$ be the $2$-handlebody obtained from $X$ by applying a $W^+$-modification $($resp.\ $W^-$-modification$)$. Then $X^+$ and $X^-$ can be embedded into $X$. Furthermore, $X$ can be embedded into $X^-$. 
\end{lemma}

Next we define Legendrian versions of $W$-modifications for Legendrian handlebodies (recall Definition~\ref{sec:stein:def:handlebody}).

\begin{figure}[ht!]
\begin{center}
\includegraphics[width=2.3in]{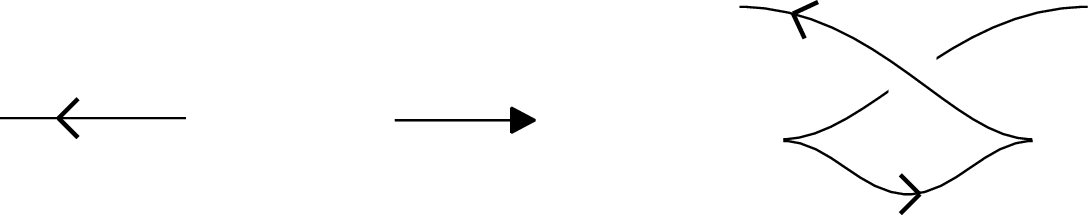}
\caption{Legendrian isotopy}
\label{fig:nuclei_fig8}
\end{center}
\end{figure}

\begin{figure}[ht!]
\begin{center}
\includegraphics[width=3.9in]{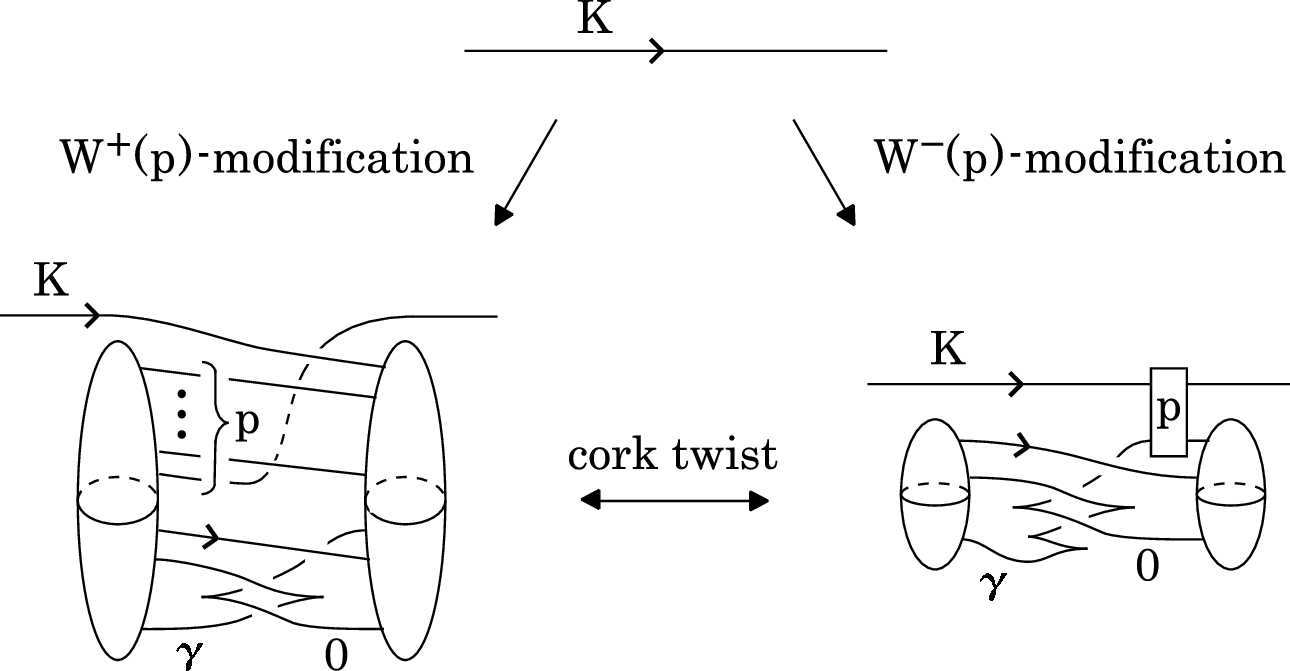}
\caption{$W^+(p)$- and $W^-(p)$-modification $(p\geq 1)$. Every framing is Seifert framing. The framing of $K$ is unchanged.}
\label{fig:nuclei_fig9}
\end{center}
\end{figure}

Let $K$ be a $2$-handle of a Legendrian handlebody. Take a small segment of the attaching circle of $K$ as in the first row of Figure~\ref{fig:nuclei_fig9}. Without loss of generality, we may assume that the orientation of the segment of $K$ is from the left to the right (Otherwise locally apply the Legendrian isotopy in Figure~\ref{fig:nuclei_fig8}. Note that this isotopy does not change the Thurston-Bennequin number and the rotation number).  

\begin{definition}\label{sec:attachment:def:attachment} Let $p$ be a positive integer. We call the local operations shown in the left and the right side of Figure~\ref{fig:nuclei_fig9} a \textit{$W^+(p)$-modification} to $K$ and a \textit{$W^-(p)$-modification} to $K$, respectively. Here we orient the $2$-handles as in the figure. Hence, each operation produces a new Legendrian handlebody from a given Legendrian handlebody. When we see Legendrian handlebodies as smooth handlebodies, these definitions and the orientations are consistent with those in Definition~\ref{sec:attachment:def:attachment} and Figure~\ref{fig:nuclei_fig7} (We can check this just by converting the $1$-handle notation). Note that the auxiliary $2$-handle $\gamma$ to any $W^+(p)$- (resp.\ $W^-(p)$-) modification satisfies the following: its framing is $0$ (resp.\ $0$); $tb(\gamma)=2$ (resp.\ $tb(\gamma)=1$); $r(\gamma)=0$ (resp.\ $r(\gamma)=1$). \end{definition}
The above definition clearly gives the following.
\begin{proposition}[\cite{AY5}]\label{sec:attachment:prop:tb and r}Let $K$ be a $2$-handle of a Legendrian handlebody. 
\\
$(1)$ Every $W^+(p)$-modification to $K$ has the following effect.
\begin{itemize}
 \item $tb(K)$ is increased by $p$, and $r(K)$ is unchanged. 
\end{itemize}
$(2)$ Every $W^-(p)$-modification to $K$ has the following effect.  
\begin{itemize}
 \item $tb(K)$ and $r(K)$ are unchanged. 
\end{itemize}
\end{proposition}



\section{Exotic smooth structures on 4-manifolds with boundary}\label{sec:stein and exotic}
In this section, applying Theorem~\ref{sec:intro:thm:original}, we show that many compact $4$-manifolds with boundary admit infinitely many distinct exotic smooth structures after $W^+$-modifications. 
Namely, we prove the following, which is a restatement of Theorem~\ref{sec:intro:thm:modified}.
\begin{theorem}\label{sec:infinite:thm:modified}
Let $X$ be a $2$-handlebody which contains a nucleus $N$ as a subhandlebody. 
Suppose that $N$ admits a Stein structure. 
Then, there exists a compact connected oriented smooth $4$-manifold $X_0$ which satisfies the following.\\
$(1)$ $X_0$ admits infinitely many distinct exotic smooth structures. \\
$(2)$ The fundamental group, the integral homology groups, the integral homology groups of the boundary, and the intersection form of $X_0$ are isomorphic to those of $X$. \\
$(3)$ $X_0$ $($resp.\ $X$$)$ can be embedded into $X$ $($resp.\ $X_0$$)$. 
\end{theorem}
\begin{proof}Since $N$ admits a Stein structure, we may assume that the subhandlebody $N$ itself is a Stein handlebody (if necessary by changing the decomposition of $N$). 
 Let $L_1,L_2,\dots,L_k$ denote all the $2$-handles of $X$ except those of $N$. By isotopy of the attaching circles of $L_1,L_2,\dots,L_k$, we may assume that $X$ is a Legendrian handlebody (cf.\ Section 4.2 of \cite{OS1}). Note that we keep the Stein handle decomposition of $N$ unchanged. We next apply $W^+(p_i)$-modification  to each $L_i$ $(1\leq i\leq k)$, where $p_i$ is a sufficiently large integer. Then we can add zig-zags to each $L_i$ and each auxiliary $2$-handle so that the Seifert framings of the resulting 2-handles are one less than  their Thurston-Bennequin numbers. The resultant handlebody (say $X_1$) is clearly a Stein handlebody and contains the nucleus $N$. Hence, we can embed $X_1$ into a closed symplectic $4$-manifold with $b_2^+>1$. Theorem~\ref{sec:intro:thm:original} thus shows that $X_1$ admits infinitely many exotic smooth structures. Let $X_0$ be the $4$-manifold obtained from $X_1$ by replacing the each $W^+(p_i)$-modification $(1\leq i\leq k)$ to $W^-(p_i)$-modification.  Namely, $X_0$ is the $4$-manifold obtained from $X_1$ by cork twists along disjoint copies of $W_1$. $X_1$ is thus homeomorphic to $X_0$. Hence the claim (1) follows. Lemmas~\ref{prop:relation of attachments} and \ref{lem:W:embedding} clearly give the claims (2) and (3). 
\end{proof}
This theorem immediately gives the following. 
\begin{corollary}\label{cor:exotic:2-handlebody:boundary sum}Let $Y$ be a $2$-handlebody, and let $N$ be a nucleus which admits a Stein structure. Put $X=Y\natural N$. Then there exists a compact connected oriented smooth $4$-manifold $X_0$ which satisfies $(1)$--$(3)$ of Theorem~\ref{sec:infinite:thm:modified}. 
\end{corollary}

Since every finitely presented group (resp.\ integral symmetric bilinear form $\mathbb{Z}^k\times \mathbb{Z}^k\to \mathbb{Z}$ $(k\geq 0)$) is realized  as the fundamental group (resp.\ intersection form) of a $2$-handlebody, we have the following corollary. Note that $(1)$ also follows from a result of Park~\cite{P2}. 
\begin{corollary}\label{sec:infinite exotic:cor}$(1)$ Let $G$ be any finitely presented group. Then there exists a compact connected oriented smooth $4$-manifold with boundary such that it admits infinitely many distinct exotic smooth structures and that its fundamental group is isomorphic to $G$.\\
$(2)$ Let $Q: \mathbb{Z}^k\times \mathbb{Z}^k\to \mathbb{Z}$ $(k\geq 0)$ be any integral symmetric bilinear form, and let $R: \mathbb{Z}^2\times \mathbb{Z}^2 \to \mathbb{Z}$ be any integral, indefinite, unimodular symmetric bilinear form. Then there exists a simply connected compact oriented smooth $4$-manifold with boundary such that it admits infinitely many distinct exotic smooth structures and that its intersection form is isomorphic to $Q\oplus R$.
\end{corollary}
\begin{remark}Akbulut and the author~\cite{AY5} recently constructed finitely many smooth structures on $4$-manifolds with boundary and with $b_2\geq 1$, using $W$-modifications. However, to the author's knowledge, there are no examples of compact $4$-manifolds with $b_2=1$ which admits infinitely many exotic smooth structures. Since Akbulut~\cite{A4}  proved that knot surgery produces at least one exotic smooth structure of the cusp neighborhood $C$, it is interesting to see whether knot surgery produces infinitely many exotic smooth structures on $C$. \end{remark}
\section{3-manifolds bounding exotic 4-manifolds}\label{sec:boundary}To state the results simply and precisely, we recall the following convention. 
For a closed oriented $3$-manifold $M$ and an oriented $4$-manifold $X$, we say ``$M$ bounds $X$'' if the oriented boundary of $X$ is diffeomorphic to $M$. 

In this section, we prove that many $3$-manifolds bound (simply connected) $4$-manifolds which admit infinitely many smooth structures,  applying Theorem~\ref{sec:intro:thm:original}. Namely, we show that many $3$-manifolds bound (simply connected) $4$-manifolds satisfying the assumption of Theorem~\ref{sec:intro:thm:original}. 

A compact connected oriented $3$-manifold $M$ is said to be \textit{Stein fillable} if there exists a compact connected Stein $4$-manifold whose boundary is orientation preserving diffeomorphic to $M$. We first prove the following theorem, which is a restatement of Theorem~\ref{intro:Stein boundary}. 
\begin{theorem}\label{thm:stein fillable exotic}$(1)$ Every Stein fillable $3$-manifold bounds a simply connected compact oriented smooth $4$-manifold which admits infinitely many distinct smooth structures. \\
$(2)$ Let $n$ be any positive integer, and let $M_1,M_2,\dots,M_n$ be any Stein fillable $3$-manifolds. Then the disjoint union $M=\coprod_{i=1}^n M_i$ bounds a simply connected compact oriented smooth $4$-manifold which admits infinitely many distinct smooth structures.
\end{theorem}
\begin{proof}$(1)$ Let $M$ be a Stein fillable $3$-manifold, and let $\widetilde{M}$ be a Stein handlebody with $\partial \widetilde{M}\cong M$. Convert the $1$-handle notation of the handlebody diagram of $\widetilde{M}$ into the dotted circle notation. This naturally gives a framed link $L$ in $S^3$ such that the Dehn surgery along $L$ gives $M$. Here recall the Wirtinger presentation of $\pi_1(S^3-\textnormal{int}\,\nu(L))$. This group is generated by a link (say $\mathbf{m}$) which consists of meridians of $L$ (cf.\ Exercise 5.2.2.(b) of~\cite{GS}).  Since $M$ is obtained from $S^3-\textnormal{int}\,\nu(L)$ by attaching $3$-dimensional $2$- and $3$-handles, the link $\mathbf{m}$ in $M$ also generates $\pi_1(M)$. Attach $4$-dimensional $2$-handles to $\widetilde{M}$ along a Legendrianization of the link $\mathbf{m}$ so that the resulting handlebody (say $Y$) is a Stein handlebody. This construction shows $\pi_1(Y-\textnormal{int}\, \widetilde{M})\cong 1$. 

Let $N$ be a nucleus which admits a Stein structure (e.g.\ Gompf nucleus $G(n)$ $(n\geq 2)$). The boundary sum $Y\natural N$ clearly admits a Stein handle decomposition. We can thus embed $Y\natural N$ into a closed symplectic $4$-manifold $Z$ with $b_2^+>1$, using a method of Akbulut-Ozbagci in~\cite{AO3}. Their method (cf.\ the proof of Theorem~\ref{sec:stein:th:closing of Stein:symplectic}) allows us to assume that $Z$ is obtained from $Y\natural N$ by attaching $2$-, $3$-, and $4$-handles only. Note that we do not assume the minimality of $Z$. Since $\pi_1((Y\natural N)-\textnormal{int}\, \widetilde{M})\cong 1$, the $4$-manifold $X:=Z-\textnormal{int}\, \widetilde{M}$ is simply connected. 

Theorem~\ref{sec:intro:thm:original} shows that $X$ (hence the reverse orientation $\overline{X}$ of $X$) has infinitely many distinct exotic smooth structures, because $X$ contains the nucleus $N$ and is embedded into $Z$. Note $\partial \overline{X}\cong M$. 

$(2)$ Let $M_i$ $(1\leq i\leq n)$ and $M$ be as in the assumption, and let $\widetilde{M}_i$ $(1\leq i\leq n)$ be a Stein handlebody with $\partial \widetilde{M}_i\cong M_i$. Similarly to $(1)$, for each $1\leq i\leq n$, we obtain a Stein handlebody $Y_i$ such that $Y_i$ contains $\widetilde{M}_i$ and that $\pi_1(Y_i-\textnormal{int}\, \widetilde{M}_i)\cong 1$. Let $\widetilde{M}$ denote the disjoint union $\coprod_{i=1}^n \widetilde{M}_i$, and let $Y$ denote the boundary sum $\natural_{i=1}^n Y_i$. Let $N$ be a nucleus which admits a Stein structure. Then the boundary sum $Y\natural N$ is clearly a Stein handlebody. Note that the complement $(Y\natural N)-\textnormal{int}\, \widetilde{M}(= N\natural_{i=1}^n (Y_i-\textnormal{int}\, \widetilde{M}_i))$ is simply connected. Similarly to $(1)$, a method of Akbulut-Ozbagci in \cite{AO3} gives a closed symplectic $4$-manifold $Z$ with $b_2^+>1$ such that $Z$ contains $Y\natural N$ and that $X:=Z-\textnormal{int}\, \widetilde{M}$ is simply connected. Since $X$ contains the nucleus $N$ and is embedded into $Z$, Theorem~\ref{sec:intro:thm:original}  shows that $X$ (hence the reverse orientation $\overline{X}$) admits infinitely many distinct exotic smooth structures. Note $\partial \overline{X}\cong M$.
\end{proof}

\begin{figure}[ht!]
\begin{center}
\includegraphics[width=4.4in]{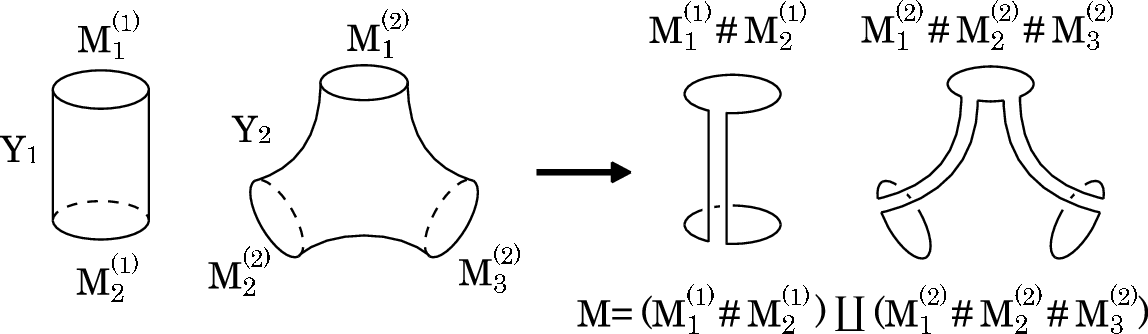}
\caption{A schematic picture of Theorem~\ref{thm:boundary:symplectic}.(2) ($\dim M=1$)}
\label{fig:nuclei_fig10}
\end{center}
\end{figure}

We extend this theorem to more general $3$-manifolds, ignoring simple connectivity of 4-manifolds. Here we briefly recall weakly symplectically fillable $3$-manifolds. For details, the reader can consult~\cite{OS1}. 
A compact symplectic manifold $(X, \omega)$ is said to have \textit{weakly convex boundary} if its boundary $\partial X$ admits a contact structure $\xi$ such that $\omega|_\xi > 0$. A closed oriented $3$-manifold $M$ is said to be \textit{weakly symplectically fillable} if $M$ bounds a compact symplectic $4$-manifold with weakly convex boundary. Such a symplectic $4$-manifold is called a \textit{weak symplectic filling} of $M$. It is known that a compact Stein $4$-manifold is a symplectic 4-manifold with weakly convex boundary. Stein fillable $3$-manifolds are thus weakly symplectically fillable. For the following Theorem~\ref{thm:boundary:symplectic}.(2),  see also the schematic picture in Figure~\ref{fig:nuclei_fig10}.
\begin{theorem}\label{thm:boundary:symplectic}
$(1)$ Every connected weakly symplectically fillable $3$-manifold bounds a compact connected oriented smooth $4$-manifold which admits infinitely many distinct smooth structures.
\\
$(2)$ Let $n$ and $m_1, m_2, \dots, m_n$ be any positive integers. Let each $M^{(i)}_j$ $($$i\in\{1,2, \dots, n\}$, \: $j\in\{1, 2, \dots m_i\}$$)$ be a connected oriented closed $3$-manifold. 
Let $M$ denote the disjoint union $\coprod_{i=1}^n (\#_{j=1}^{m_i} M^{(i)}_j)$. Assume that, for each $1\leq i\leq n$, the disjoint union $\coprod_{j=1}^{m_i} M^{(i)}_j$ bounds a compact connected oriented smooth $4$-manifold $Y_i$. Assume further that the disjoint union $Y:=\coprod_{i=1}^{n} Y_i$ can be embedded into a compact connected symplectic $4$-manifold with weakly convex boundary. Then $M$ bounds a compact connected oriented smooth $4$-manifold which admits infinitely many distinct smooth structures. 
\end{theorem}
\begin{proof}
$(1)$ Let $M$ be a weakly symplectically fillable $3$-manifold, and let $\widetilde{M}$ be a weak symplectic filling of $M$. 
The contact $3$-manifold $\partial \widetilde{M}=M$ is represented by a contact $(\pm 1)$-surgery along a Legendrian link $L$ in the standard contact $S^3$ (\cite{DG}). Attach two $2$-handles to $\partial \widetilde{M}$ along the framed Legendrian link as shown in Figure~\ref{fig:nuclei_fig4}, where we draw the link in a $3$-ball away from $L$ in $S^3$ (the framings in the figure are the Seifert framings.). The resulting $4$-manifold is clearly the boundary sum $\widetilde{M}\natural G(n)$ $(n\geq 2)$. According to a theorem of Weinstein \cite{Wei} (see also Theorem~5.8 of \cite{Et2}), this construction ensures that $\widetilde{M}\natural G(n)$ admits a symplectic structure with weakly convex boundary. By a theorem of Eliashberg \cite{E3} and Etnyre \cite{Et}, we can embed $\widetilde{M}\natural G(n)$ into a closed symplectic $4$-manifold $Z$ with $b_2^+>1$. (The $b_2^+>1$ condition is easily verified by using an embedding of $\widetilde{M}\natural G(n)\natural G(n)$.) Theorem~\ref{sec:intro:thm:original} thus shows that $X:=Z-\textnormal{int}\, \widetilde{M}$ admits infinitely many distinct exotic smooth structures. Note $\partial \overline{X}\cong M$.

$(2)$  Let $S$ be a compact connected symplectic $4$-manifold with weakly convex boundary which contains $Y$. Similarly to the proof of $(1)$, we obtain a closed connected symplectic $4$-manifold $Z$ with $b_2^+>1$ which contains $S\natural G(k)$ $(k\geq 2)$. 
Put $M_i=\#_{j=1}^{m_i} M^{(i)}_j$ $(1\leq i\leq n)$ and $Z^0=Z-\textnormal{int}\, Y$. 
Since each $Y_i$ $(1\leq i\leq n)$ is path connected, we can take a path $\gamma^{(i)}_j\subset Y_i$ which connects points of $M^{(i)}_{1}$ and $M^{(i)}_j$, for each $2\leq j\leq m_i$. By dimensional reason, we may assume that $\gamma^{(i)}_2, \gamma^{(i)}_3, \dots, \gamma^{(i)}_{m_i}$ are smoothly embedded mutually disjoint paths in $Y_i$. 
Taking their regular neighborhoods, we get the corresponding $(m_i-1)$ mutually disjoint $4$-dimensional $1$-handles $\widehat{\gamma}^{(i)}_2, \widehat{\gamma}^{(i)}_3, \dots, \widehat{\gamma}^{(i)}_{m_i}\subset Y_i$ for each $1\leq i\leq n$. The attaching region of each $\widehat{\gamma}^{(i)}_j$ is located in $M^{(i)}_{1}\coprod M^{(i)}_j$. Note that these 1-handles naturally give the connected sum $M_i=\#_{j=1}^{m_i} M^{(i)}_j$ embedded in $Y_i$. Using these $1$-handles $\widehat{\gamma}^{(i)}_2, \widehat{\gamma}^{(i)}_3, \dots, \widehat{\gamma}^{(i)}_{m_i}$ $(1\leq i\leq n)$, we join connected components of $Z^0$ in $Z$. Let $X$ be the resulting compact oriented smooth $4$-manifold. This construction shows that, for each $1\leq i\leq n$, $M^{(i)}_1, M^{(i)}_2,\dots, M^{(i)}_{m_i}$ belong to the same connected component of $X$. Since $\partial Y_i=\coprod_{j=1}^{m_i} M^{(i)}_j$, this shows that the number of connected components of $X$ is the same as that of $X\cup (\coprod _{i=1}^{n} Y_i)$. Thus $X$ is connected, because $X\cup (\coprod _{i=1}^{n} Y_i)=Z$. Note that $X$ contains the nucleus $N$ and that $\partial X=\overline{M}$. Theorem~\ref{sec:intro:thm:original} thus shows that $X$ admits infinitely many distinct smooth structures. 
\end{proof}
\begin{corollary}\label{cor:symplectic:boundary}Let $M$ be a weakly symplectically fillable 3-manifold. Then $M\#\overline{M}$ bounds a compact connected oriented smooth $4$-manifold which admits infinitely many distinct smooth structures. 
\end{corollary}
\begin{proof}Let $\widetilde{M}$ be a weak symplectic filling of $M$. Note that $Y=M\times [0,\,1]$ is embedded into $\widetilde{M}$ and that $\partial Y=M\coprod \overline{M}$. Thus, the desired claim follows from Theorem~\ref{thm:boundary:symplectic}.(2). 
\end{proof}
This corollary provides non weakly symplectically fillable 3-manifolds bounding $4$-manifolds which admit infinitely many smooth structures. 
\begin{example}\label{ex:symplectic:boundary} Let $\Sigma(2,3,5)$ denote the Poincar\'{e} homology sphere. Since $\Sigma(2,3,5)$ is Stein fillable, the above corollary shows that $\Sigma(2,3,5)\#\overline{\Sigma(2,3,5)}$ bounds  a compact connected oriented smooth $4$-manifold which admits infinitely many distinct smooth structures. Note that $\Sigma(2,3,5)\#\overline{\Sigma(2,3,5)}$ is not weakly symplectically fillable (hence, non Stein fillable) for both orientations (\cite{L}).
\end{example}

\section{Non-existence of Stein structures}\label{sec:non-existence}
In this section, we show that log transform and knot surgery do not produce compact Stein $4$-manifolds under a mild condition. 

\begin{figure}[ht!]
\begin{center}
\includegraphics[width=2.0in]{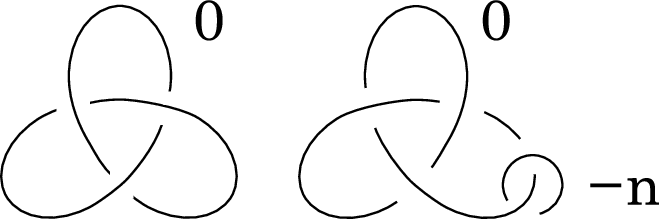}
\caption{$\overline{C}\natural G(n)$}
\label{fig:nuclei_fig11}
\end{center}
\end{figure}

\begin{lemma}\label{nonstein:lem:nonzeroclass}Let $X$ be a compact connected oriented smooth $4$-manifold $($possibly with boundary$)$ which contains a torus $T$. Let $X_{(p)}$ $(p\geq 1)$ and $X_K$ be the $p$-log transform along $T$ and the knot surgery along $T$ with a knot $K$, respectively. Then the following hold. \\
$(1)$ Assume that $T$ is a $c$-embedded torus. Then $[T]$ is non-torsion in $H_2(X;\mathbb{Z})$ if and only if $[T_p]$ is non-torsion in $H_2(X_{(p)};\mathbb{Z})$. \\
$(2)$ $[T]$ is non-torsion in $H_2(X;\mathbb{Z})$ if and only if $[T]$ is non-torsion in $H_2(X_K;\mathbb{Z})$, where the last $[T]$ denotes the class of a parallel copy of the original $T$. $($This holds even when $T\subset X$ is not $c$-embedded.$)$

\end{lemma}
\begin{proof} $(1)$ Applying the $p$-log transform procedure in Figure~17 of \cite{AY1} to the cusp neighborhood $C$, we easily see $\pi_1(C_{(p)})\cong \mathbb{Z}/p\mathbb{Z}$. The Mayer-Vietoris exact sequences for $X=(X-\textnormal{int}\, C)\cup C$ and $X_{(p)}=(X-\textnormal{int}\, C_{(p)})\cup C_{(p)}$ thus imply $b_1(X_{(p)})=b_1(X)$. 

Suppose that $[T]$ is non-torsion in $H_2(X;\mathbb{Z})$. We show that $[T_p]$ is non-torsion in $H_2(X_{(p)};\mathbb{Z})$. 
It suffices to see that the homomorphism $H_2(\nu (T_p);\mathbb{Q})\to H_2(X_{(p)};\mathbb{Q})$ induced by the inclusion is injective.

(i) The case where $X$ is closed. Suppose that the above homomorphism is not injective, thus, a 0-map. Then we can easily show $b_3(X_{(p)})=b_3(X_{(p)})-1$ by using two homology exact sequences for the pairs $(X, \nu (T))$ and $(X_{(p)}, \nu (T_p))$ (Note $H_*(X, \nu (T))\cong H_*(X_{(p)}, \nu (T_p))$). Since $X$ and $X_{(p)}$ are closed 4-manifolds, we have $b_1(X_{(p)})=b_3(X_{(p)})$ and $b_1(X)=b_3(X)$. We thus get $b_1(X_{(p)})\neq b_1(X)$, which is a contradiction. Hence the homomorphism $H_2(\nu (T_p);\mathbb{Q})\to H_2(X_{(p)};\mathbb{Q})$ is injective. 

(ii) The case where $X$ has non-empty boundary. Take a double $DX$ of $X$ to form a closed $4$-manifold. The Mayer-Vietoris exact sequence shows that the inclusion induced homomorphism $H_2(X;\mathbb{Q})\to H_2(DX;\mathbb{Q})$ is injective. The required claim in the closed case thus shows that the inclusion induced homomorphism $H_2(\nu (T_p);\mathbb{Q})\to H_2((DX)_{(p)};\mathbb{Q})$ is injective. Therefore, the homomorphism $H_2(\nu (T_p);\mathbb{Q})\to H_2(X_{(p)};\mathbb{Q})$ is injective. 

Conversely, suppose that $[T_p]$ is non-torsion in $H_2(X_{(p)};\mathbb{Z})$. Then the same argument as above shows that $H_2(\nu (T);\mathbb{Q})\to H_2(X;\mathbb{Q})$ is injective. 

(2) By the gluing condition of knot surgery, we can easily check $b_1(X_K)=b_1(X)$ using the Mayer-Vietoris exact sequences. Let $T_K$ be the torus given in Definition~\ref{def:knot surgery}. The natural inclusion induces an isomorphism $H_2(\nu (T_K);\mathbb{Q})\to H_2((S^3-\nu(K))\times S^1;\mathbb{Q})$. Note that $[T_K]=[T]$ in $H_2(X_K;\mathbb{Q})$ and that $H_*((S^3-\nu(K))\times S^1;\mathbb{Q})\cong H_*(T^2\times D^2;\mathbb{Q})$. The following conditions (a) and (b) are thus equivalent: (a) The inclusion induced homomorphism $H_2(\nu (T);\mathbb{Q})\to H_2(X_K;\mathbb{Q})$ is injective; (b) The inclusion induced homomorphism $H_2((S^3-\nu(K))\times S^1;\mathbb{Q})\to H_2(X_K;\mathbb{Q})$ is injective. Since $H_*(X, \nu (T))\cong H_*(X_K, (S^3-\nu(K))\times S^1)$, the rest of the proof is similar to the (1) case. 
\end{proof}
\begin{remark} If $T$ is not $c$-embedded, Lemma~\ref{nonstein:lem:nonzeroclass}.(1) does not always hold. For example, put $T=T^2\times \{pt.\}$ in $T^2\times S^2$. Then, applying the procedure in \cite{AY1}, we easily see that the $p$-log transform $(p\geq 1)$ of $T^2\times S^2$ is diffeomorphic to $S^1\times S^3$. Consequently $[T_p]=0$, while $[T]$ is non-torsion in $H_2({T^2\times S^2};\mathbb{Z})$. 
\end{remark}
\begin{lemma}\label{lem:nonstein:reverse}
Let $X$ be a compact connected oriented smooth $4$-manifold with boundary which contains a cusp neighborhood. Then the reverse orientation $\overline{X}$ of $X$ does not admit any Stein structure. 
\end{lemma}
\begin{proof}
By the assumption, the boundary sum $\overline{X}\natural G(n)$ $(n\geq 2)$ contains the boundary sum $\overline{C}\natural G(n)$, where $G(n)$ and $\overline{C}$ denote the Gompf nucleus and the reverse orientation of $C$, respectively. A handlebody picture of $\overline{C}\natural G(n)$ is given in Figure~\ref{fig:nuclei_fig11}. Sliding the right trefoil knot over the left trefoil knot, we get a $0$-framed ribbon knot. This ribbon knot gives a smoothly embedded sphere $\Sigma\subset \overline{C}\natural G(n)$ with $\Sigma\cdot \Sigma=0$. 
  Now suppose that $\overline{X}$ admits a Stein structure. Then $\overline{X}\natural G(n)$ admits a Stein structure, because $G(n)$ also admits a Stein structure (Example~\ref{ex:gompf:stein}). We can thus embed the Stein $4$-manifold $\overline{X}\natural G(n)$ into a closed symplectic $4$-manifold $Z$ with $b_2^+>1$. Note that  the sphere $\Sigma$ represents a non-torsion class of $H_2(Z;\mathbb{Z})$, because it algebraically intersects with the obvious sphere with the self-intersection number $-n$. Since $[\Sigma]\cdot [\Sigma]=0$, this contradicts the adjunction inequality. 
\end{proof}

The following is a restatement of Theorem~\ref{intro:non-stein}. 
\begin{theorem}\label{sec:nonstein:thm:nonstein}Let $X$ be a compact connected oriented smooth $4$-manifold with boundary which contains a $c$-embedded torus $T$. Assume that $[T]$ is non-torsion in $H_2(X;\mathbb{Z})$. Let $X_{(p)}$ and $X_K$ denote the $p$-log transform along $T$ and the knot surgery along $T$ with a knot $K$, respectively. Then the following hold.\\
$(1)$ For $p\geq 2$, $X_{(p)}$ does not admit any Stein structure for both orientations.\\
$(2)$ For any knot $K$ in $S^3$ with the non-trivial Alexander polynomial, $X_K$ does not admit any Stein structure for both orientations. 
\end{theorem}
We give two proofs of this theorem, using two different embeddings of a Stein $4$-manifold. Namely, Lisca-Mati\'{c}'s embedding into a minimal complex surfaces of general type and Akbulut-Ozbagci's embedding into a closed symplectic $4$-manifold. While our second proof is a little more technical, we hope that the second proof is useful to show non-existence of some other structures under some conditions, since a symplectic $4$-manifold with weakly convex boundary also can be embedded into a closed symplectic $4$-manifold (\cite{E3}, \cite{Et}). 
\begin{proof}[Proof 1 of Theorem~\ref{sec:nonstein:thm:nonstein}]$(1)$ By Lemma~\ref{nonstein:lem:nonzeroclass}, $[T_p]$ is non-torsion in $H_2(X_{(p)};\mathbb{Z})$. Suppose that $X_{(p)}$ with the orientation induced from $X$ admits a Stein structure. Then $X_{(p)}$ has a Stein handlebody decomposition. Attach a $2$-handle along the meridian of the attaching circle of each $2$-handle of the Stein handlebody $X_{(p)}$ so that the resulting manifold (say $Y_{(p)}$) is a Stein handlebody. 
Theorem~\ref{sec:stein:th:closing of Stein:complex} thus gives a minimal complex surface (say $Z_{(p)}$) of general type and with $b_2^+>1$ which contains $Y_{(p)}\supset X_{(p)}$. The construction of $Y_{(p)}$ ensures that $[T_p]$ is non-torsion in $H_2(Z_{(p)};\mathbb{Z})$. 

Since $X_{(p)}$ is the $p$-log transform of $X$ along $T$, reversing this operation along $T_p$ in $Z_{(p)}$, we obtain a closed oriented smooth 4-manifold (say $Z$) which contains $X$. Note that $Z_{(p)}$ is the $p$-log transform of $Z$ along the $c$-embedded torus $T$ in $X\subset Z$. 
Since $[T_p]$ is non-torsion in $H_2(Z_{(p)};\mathbb{Z})$, Lemma~\ref{nonstein:lem:nonzeroclass} shows that $[T]$ is non-torsion in $H_2(Z;\mathbb{Z})$. It is easy to see $b_2^+(Z)>1$. Theorem~\ref{Thm:SW:log} thus gives a class $K$ of $H^2(Z_{(p)};\mathbb{Z})$ such that 
\begin{equation*}
L_1:=K+(p-1)PD([T_p])\quad \text{and}\quad L_2:=K-(p-1)PD([T_p])
\end{equation*}
 are Seiberg-Witten basic classes of $Z_{(p)}$, because the Seiberg-Witten invariant of the complex surface $Z_{(p)}$ does not vanish. Since $[T_p]$ is non-torsion, we get $L_1\neq L_2$. 
 Since the minimal complex surface $Z$ of general type has only one basic class up to sign (Remark~\ref{rem:SW of general type}), we have $L_1 + L_2 = 0$. This shows $2K=0$. Applying the adjunction inequality to $L_1$ and $[T_p]$, we get $\langle K,\, [T_p] \rangle=0$. The square of $L_1$ is thus zero. This is a contradiction (see Remark~\ref{rem:SW of general type}). Hence $X_{(p)}$ does not admit any Stein structure. Since $X_{(p)}$ still contains (a smaller copy of) the cusp neighborhood, Lemma~\ref{lem:nonstein:reverse} shows that the reverse orientation $\overline{X}_{(p)}$ of $X_{(p)}$ admits no Stein structures. 
    
 The proof of (2) is similar to the above. 
\end{proof}

\begin{figure}[ht!]
\begin{center}
\includegraphics[width=3.2in]{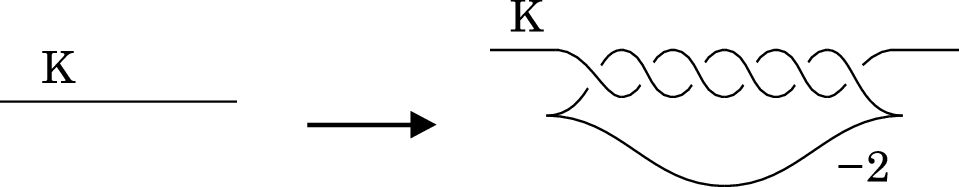}
\caption{}
\label{fig:nuclei_fig12}
\end{center}
\end{figure}

\begin{proof}[Proof 2 of Theorem~\ref{sec:nonstein:thm:nonstein}]$(1)$ By Lemma~\ref{nonstein:lem:nonzeroclass}, $[T_p]$ is non-torsion in $H_2(X_{(p)};\mathbb{Z})$. Suppose that $X_{(p)}$ with the orientation induced from $X$ admits a Stein structure. Then $X_{(p)}$ has a Stein handlebody decomposition. For each $2$-handle $K$ of the Stein handlebody $X_{(p)}$, take a small segment of the attaching circle and attach a $2$-handle as in Figure~\ref{fig:nuclei_fig12}. The resulting manifold (say $Y_{(p)}$) is clearly a Stein handlebody. Theorem~\ref{sec:stein:th:closing of Stein:symplectic} thus gives a closed symplectic $4$-manifold (say $Z_{(p)}$) with $b_2^+>1$ which contains $Y_{(p)}$. The construction of $Y_{(p)}$ ensures the existence of a smoothly embedded sphere $S$ in $Z_{(p)}$ which satisfies $[S]\cdot [S]=-2$ and $[S]\cdot [T_p]\geq 3$. 

Reversing the $p$-log transform operation along $T_p$ in $Z_{(p)}$, we obtain a closed oriented smooth 4-manifold (say $Z$) which contains $X$. Note that $Z_{(p)}$ is the $p$-log transform of $Z$ along the $c$-embedded torus $T$ in $X\subset Z$. 
 Since the Seiberg-Witten invariant of the symplectic $4$-manifold $Z_{(p)}$ does not vanish, Theorem~\ref{Thm:SW:log} gives a class $K$ of $H^2(Z_{(p)};\mathbb{Z})$ such that $K+(p-1)PD([T_p])$ and $K-(p-1)PD([T_p])$ are Seiberg-Witten basic classes of $Z_{(p)}$. Applying the adjunction inequality to these basic classes and the sphere $S$, we get the following inequalities:
\begin{align*}
 &\left|\, \langle K,\, [S] \rangle + (p-1)[T_p]\cdot S \, \right|-2\leq 0, \\
 &\left|\, \langle K,\, [S] \rangle - (p-1)[T_p]\cdot S \, \right|-2\leq 0. 
\end{align*}
This gives $-2\leq [S]\cdot [T_p] \leq 2$, which contradicts the fact $[S]\cdot [T_p]\geq 3$. Hence $X_{(p)}$ does not admit any Stein structure. Since $X_{(p)}$ still contains (a smaller copy of) the cusp neighborhood, Lemma~\ref{lem:nonstein:reverse} shows that the reverse orientation $\overline{X}_{(p)}$ of $X_{(p)}$ admits no Stein structures. 
 
 The proof of (2) is similar to (1). 
\end{proof}
\begin{remark}$(1)$ In the $p\geq 3$ (resp.\ $\deg(\triangle _K)\geq 2$) case, we can simplify the first proof  as follows. Since $[T_p]$ (resp.\ $[T]$) is non-torsion in $H_2(Z_{(p)};\mathbb{Z})$ (resp.\ $H_2(Z_{K};\mathbb{Z})$), Theorem~\ref{Thm:SW:log} (resp.\ Theorem~\ref{Thm:SW:knot}) implies that the number of Seiberg-Witten basic classes of $Z_{(p)}$ (resp.\ $Z_K$) is larger than two. This is a contradiction. \\
$(2)$ If Theorem~\ref{Thm:SW:knot} holds for non $c$-embedded torus, the above proofs show that Theorem~\ref{sec:nonstein:thm:nonstein}.(2) holds even when $T$ is not $c$-embedded. On the other hand, Theorem~\ref{sec:nonstein:thm:nonstein}.(1) does not always hold when $T$ is not $c$-embedded. Actually, any log transform of $T^2\times D^2$ along the obvious torus produces $T^2\times D^2$ itself, because any self-diffeomorphism of the boundary $S^1\times S^1\times S^1$ extends to a self-diffeomorphism of its color neighborhood $S^1\times S^1\times S^1\times [0,\,1]$. $T^2\times D^2$ is known to admit a Stein structure.  
\end{remark}
\section{Stein 4-manifolds and non-Stein 4-manifolds}\label{sec:stein and non-stein}
In this section we construct arbitrary many Stein $4$-manifolds and infinitely many non-Stein $4$-manifolds which are all homeomorphic but mutually non-diffeomorphic. 
Namely, we prove the following theorem, which is a restatement of Theorem~\ref{intro:stein and non-stein}. 
\begin{theorem}\label{sec:thm:stein and non-stein}
Let $X$ be a $2$-handlebody which contains a nucleus $N$ as a subhandlebody. 
Suppose that $N$ admits a Stein structure.   
Then for each integer $n\geq 1$, there exist infinitely many compact connected oriented smooth $4$-manifolds $X_i$ $(i=0,1,2,\dots)$ which satisfies the following.\\
$(1)$ $X_i$ $(i\geq 0)$ are all homeomorphic but mutually non-diffeomorphic.\\
$(2)$ Each $X_i$ $(1\leq i\leq n)$ admits a Stein structure, and any $X_i$ $(i\geq n+1)$ admits no Stein structure.\\
$(3)$ The fundamental group, the integral homology groups, the integral homology groups of the boundary, and the intersection form of each $X_i$ $(i\geq 0)$ are isomorphic to those of $X$.\\
$(4)$ $X$ can be embedded into $X_0$.\\
$(5)$ Each $X_i$ $(0\leq i\leq n)$ can be embedded into $X$. 
\end{theorem}

Note that Akbulut and the author~\cite{AY5} recently constructed arbitrary many Stein $4$-manifolds and arbitrary many non-Stein $4$-manifolds which are all homeomorphic but mutually non-diffeomorphic. (The following claim is a simplification of the claims of Theorem~10.1 and Remark~10.7 of \cite{AY5}. While $(4)$ and (5) are not clearly stated in~\cite{AY5}, we can easily check from the proof of Theorem~10.1 in \cite{AY5}.)
\begin{theorem}[Theorem~10.1 and Remark~10.7 of \cite{AY5}]
Let $X$ be a $2$-handlebody with $b_2\geq 1$, and let $Y$ be the boundary sum $X\natural (S^2\times D^2)$. Then, for each $n\geq 1$, there exist compact connected oriented smooth $4$-manifolds $Y_0, Y_1, \dots, Y_{2n}$ which satisfies the following. $($Furthermore, the same claim also holds in the $Y=X\natural (\overline{\mathbb{C}\mathbb{P}^2}-\textnormal{int}\, D^4)$ case.$)$\\
$(1)$ $Y_i$ $(0\leq i\leq 2n)$ are mutually homeomorphic but non-diffeomorphic.\\
$(2)$ Every $Y_i$ $(1\leq i\leq n)$ admits a Stein structure, and any $Y_j$ $(n+1\leq j\leq 2n)$ admits no Stein structure.\\
$(3)$ The fundamental groups, the integral homology groups, the integral homology groups of the boundary, and the intersection forms of every $Y_i$ $(0\leq i\leq 2n)$ are isomorphic to those of $Y$.\\
$(4)$ $Y$ can be embedded into $Y_{0}$.\\
$(5)$ Each $Y_i$ $(0\leq i\leq 2n)$ can be embedded into $Y$.
\end{theorem}

To prove Theorem~\ref{sec:thm:stein and non-stein}, we recall the following theorem of Akbulut and the author~\cite{AY5}. 
\begin{theorem}[Theorem~1.1 of \cite{AY5}]\label{thm:AY5}Let $X$ be any $2$-handlebody with $b_2\geq 1$. Then for each integer $n\geq 1$, there exist compact connected oriented smooth $4$-manifolds $X_0, X_1, \dots, X_n$ which satisfies the following.\\ 
$(1)$ The fundamental group, the integral homology groups, the integral homology groups of the boundary, and the intersection form of each $X_i$ $(0\leq i\leq n)$ are isomorphic to those of $X$.\\
$(2)$ $X_i$ $(0\leq i\leq n)$ are all homeomorphic but mutually non-diffeomorphic.\\
$(3)$ Each $X_i$ $(1\leq i\leq n)$ admits a Stein structure.\\
$(4)$ $X$ can be embedded into $X_{0}$.\\
$(5)$ Each $X_i$ $(0\leq i\leq n)$ can be embedded into $X$. 
\end{theorem}
We here briefly summarize the construction of the above $X_i$'s in Section~5 of \cite{AY5}, emphasizing not their Stein handlebody decompositions but their smooth handlebody decompositions. 
Recall the following terminology of \cite{AY5}.
\begin{definition}\label{def:goodstein}
$(1)$ A Legendrian handlebody $X$ is said to be a good Legendrian handlebody if the following two conditions are satisfied. 
\begin{enumerate}
 \item [(i)] Some $2$-handles $K_1,K_2,\dots,K_s$ of the handlebody $X$ does not algebraically go over any $1$-handles.
 \item [(ii)] The homology classes given by $2$-handles $K_1,K_2,\dots,K_s$ span a basis of $H_2(X;\mathbb{Z})$. 
\end{enumerate}
$(2)$ A Stein handlebody is said to be a good Stein handlebody if it satisfies the above conditions (i) and (ii). 
\end{definition}

\begin{figure}[ht!]
\begin{center}
\includegraphics[width=3.6in]{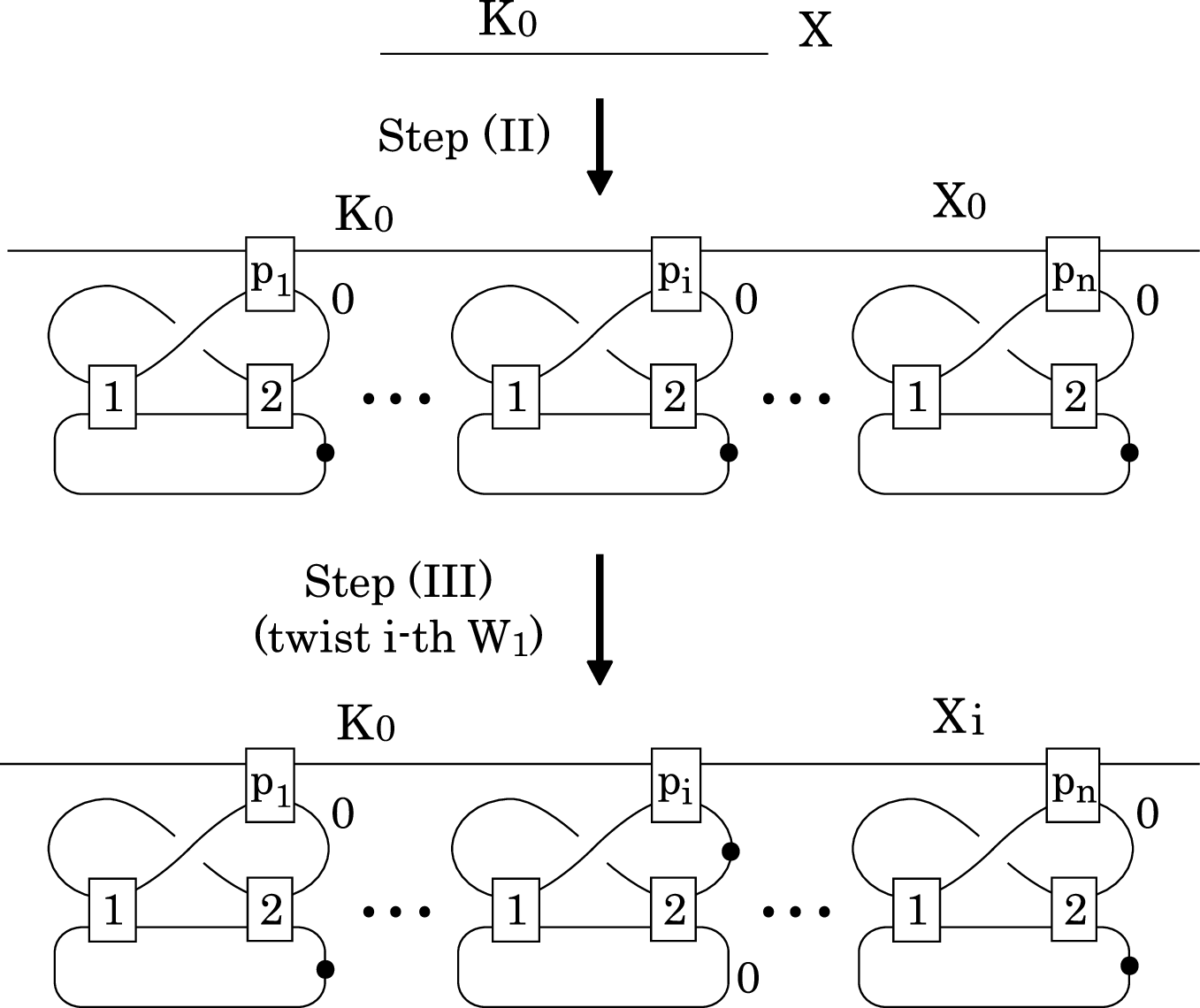}
\caption{Steps (II) and (III) to $K_0$ (ignoring Legendrian diagrams)}
\label{fig:nuclei_fig13}
\end{center}
\end{figure}

\begin{figure}[ht!]
\begin{center}
\includegraphics[width=4.5in]{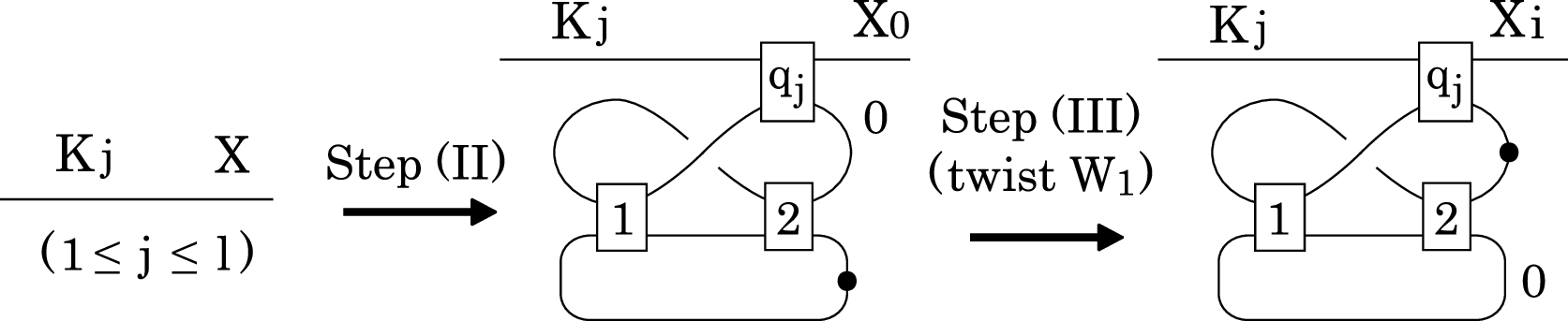}
\caption{Steps (II) and (III) to $K_j$ (ignoring Legendrian diagrams)}
\label{fig:nuclei_fig14}
\end{center}
\end{figure}

Each $X_i$ $(0\leq i\leq n)$ in Theorem~\ref{thm:AY5} is obtained from $X$, roughly, by the following Steps (I)--(III). (These steps are simplified and little different from the original construction. However, the resulting Legendrian handlebodies $X_0,X_1,\dots,X_n$ are the same.)
\begin{enumerate}
 \item [(I)] By sliding $2$-handles and isotopy, change the given $2$-handlebody $X$ into a good Legendrian handlebody.  Let $K_0,K_1,\dots,K_l$ denote all the $2$-handles of $X$. We assume that $K_0$ does not go over any $1$-handles. 
 \item [(II)] Apply $W^-(p_1)$-, $W^-(p_2)$-, \dots, $W^-(p_n)$-modifications to $K_0$ (see Figure~\ref{fig:nuclei_fig13}). Here we may choose $p_1$ to be an arbitrary large integer. Also, we may choose each $p_{i}$ $(2\leq i\leq n)$ so that $p_{i}-p_{i-1}$ is arbitrary large. Next apply a $W^-(q_j)$-modification to each $2$-handle $K_j$ $(1\leq j\leq l)$ (see Figure~\ref{fig:nuclei_fig14}). 
 We may choose each $q_j$ $(1\leq j\leq l)$ to be an arbitrary large integer. Let $X_0$ denote the resulting Legendrian handlebody. (Beware that this $X_0$ corresponds to $X^{(n)}_{-1}$ in Section~5 of \cite{AY5}.)
 \item [(III)] Fix an integer $i$ with $1\leq i\leq n$. Replace the $W^-(p_i)$-modification applied to $K_0$ with the corresponding $W^+(p_i)$-modification (i.e.\ twist the cork $W_1$, see Figure~\ref{fig:nuclei_fig13}). Also, for each $1\leq j\leq l$, replace the $W^-(q_j)$-modification applied to $K_j$ with the corresponding $W^+(q_j)$-modification (i.e.\ twist the cork $W_1$, see Figure~\ref{fig:nuclei_fig14}.). Finally add zig-zags to all the $2$-handles so that the handlebody becomes a Stein handlebody and that it satisfies certain conditions with respect to the Thurston-Bennequin numbers and the rotation numbers. The result of $X$ is the $X_i$ in Theorem~\ref{thm:AY5}. 
\end{enumerate}

Keeping these steps in mind, we prove Theorem~\ref{sec:thm:stein and non-stein}. 

\begin{figure}[ht!]
\begin{center}
\includegraphics[width=4.7in]{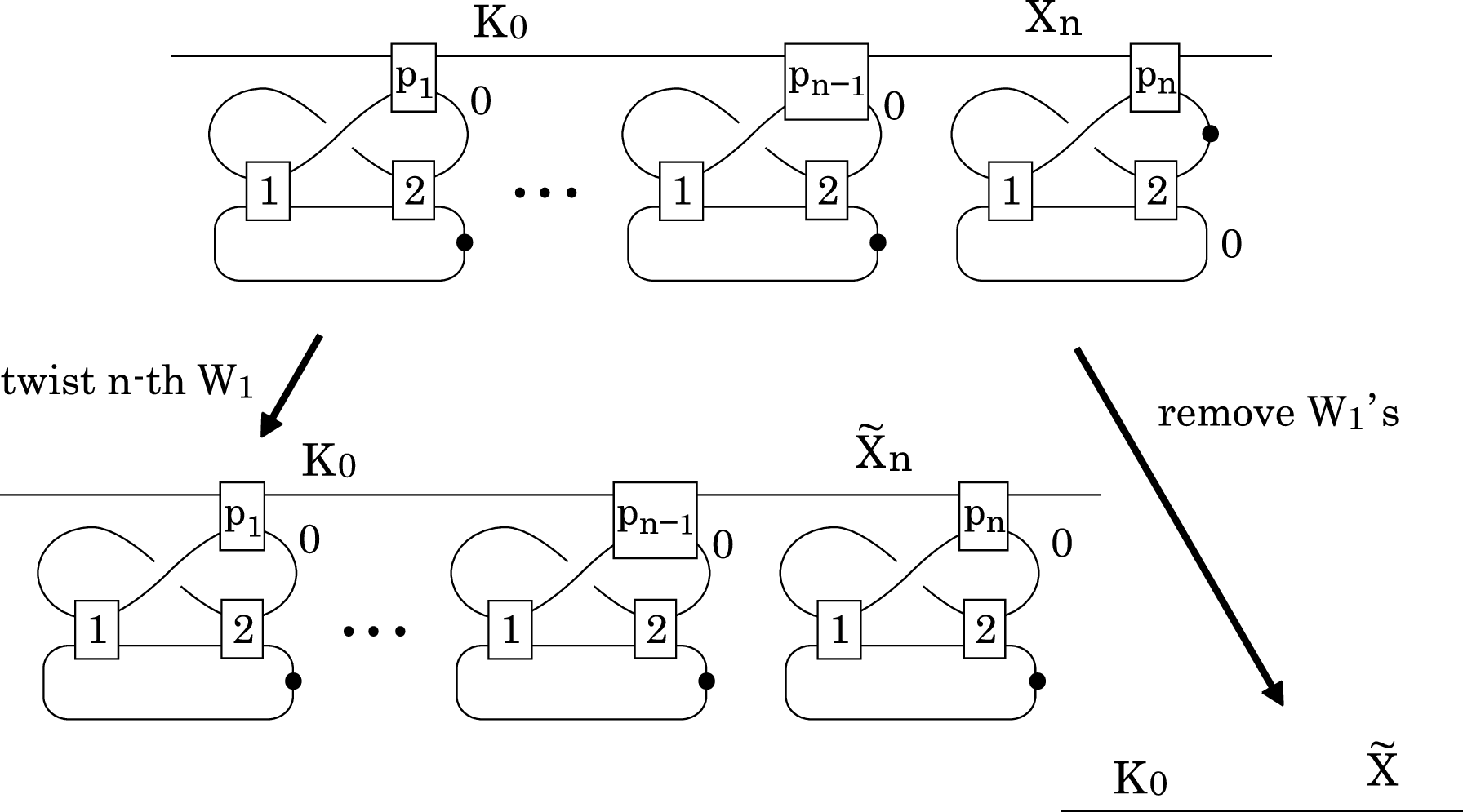}
\caption{Constructions of $\widetilde{X}_n$ and $\widetilde{X}$ (ignoring Legendrian diagrams)}
\label{fig:nuclei_fig15}
\end{center}
\end{figure}

\begin{figure}[ht!]
\begin{center}
\includegraphics[width=3.4in]{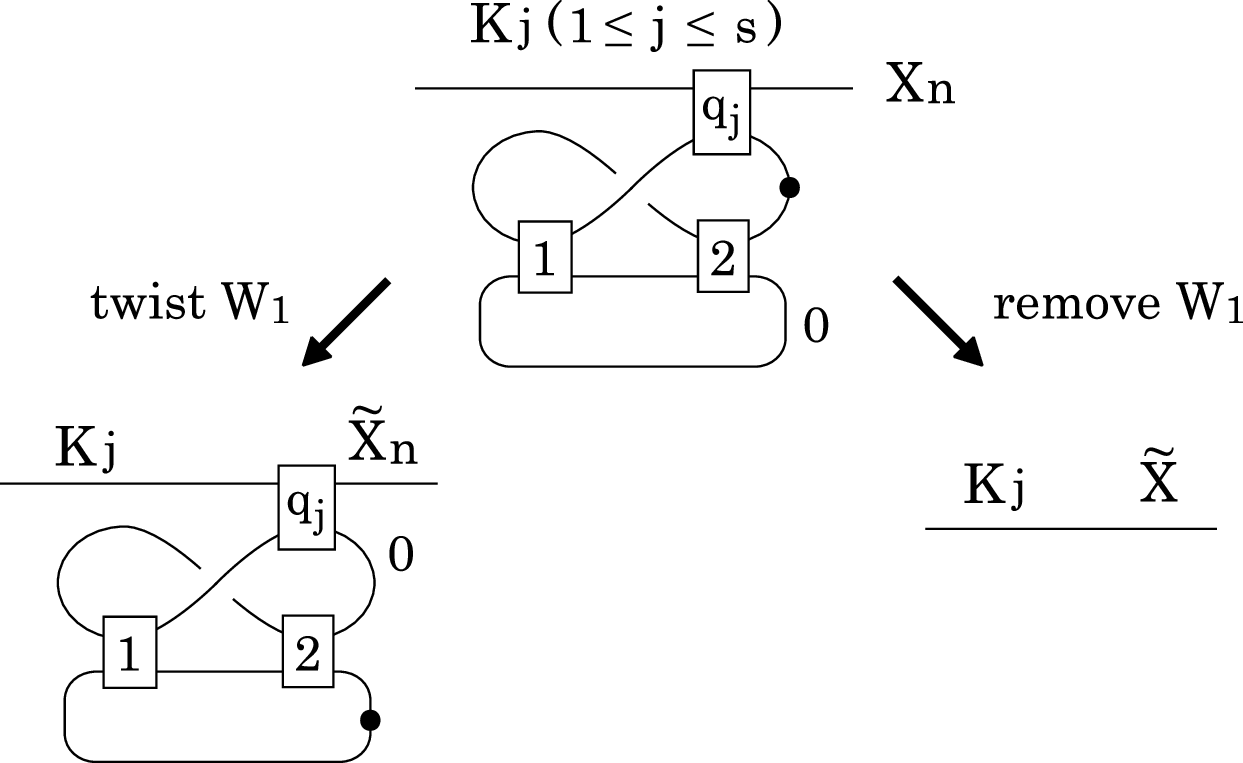}
\caption{Constructions of $\widetilde{X}_n$ and $\widetilde{X}$ (ignoring Legendrian diagrams)}
\label{fig:nuclei_fig16}
\end{center}
\end{figure}

\begin{figure}[ht!]
\begin{center}
\includegraphics[width=3.4in]{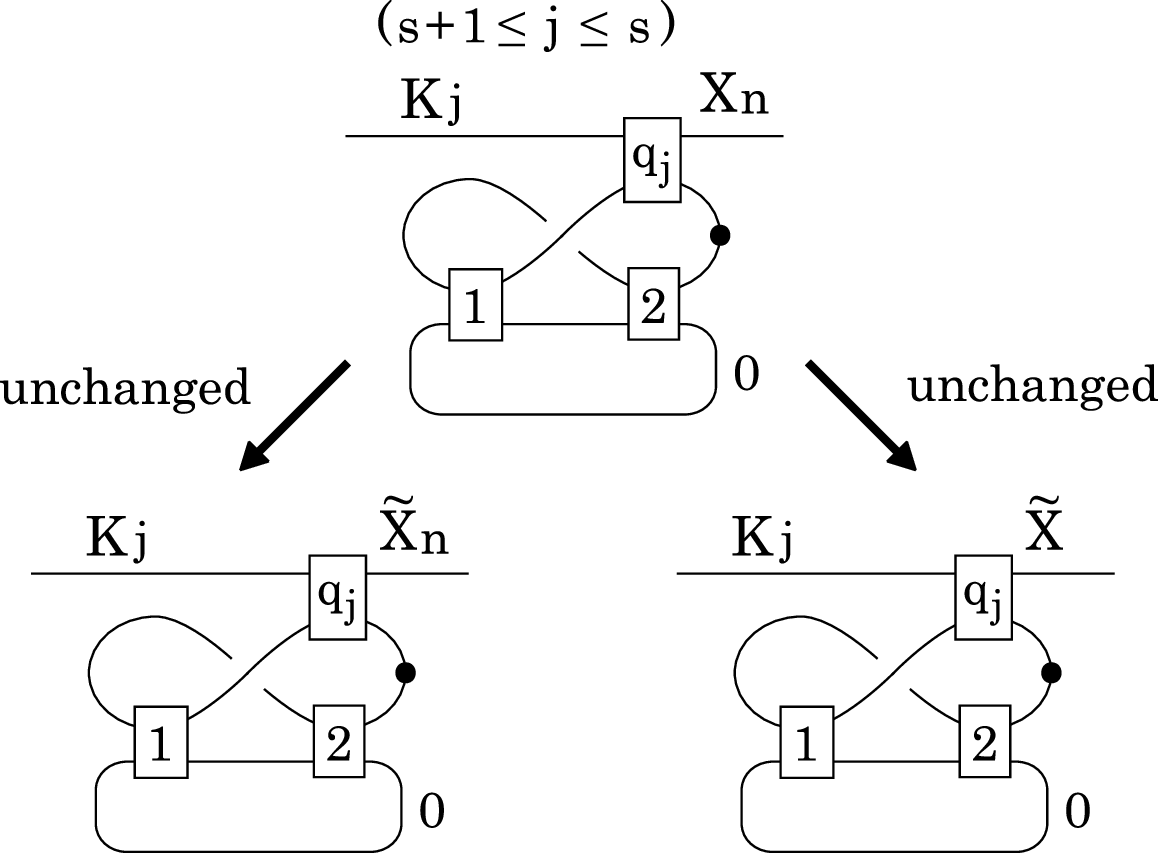}
\caption{Constructions of $\widetilde{X}_n$ and $\widetilde{X}$ (ignoring Legendrian diagrams)}
\label{fig:nuclei_fig17}
\end{center}
\end{figure}

\begin{proof}[Proof of Theorem~\ref{sec:thm:stein and non-stein}]Let $X, N, n$ be as in the assumption. Since $N$ admits a Stein structure, we may assume that the subhandlebody $N$ itself is a Stein handlebody (if necessary by changing the decomposition of $N$). 
 We carefully apply Step (I) to $X$ as follows. Let $K_0,K_1,\dots,K_s$ denote all the $2$-handles of the Stein handlebody $N$, and let $K_{s+1}, K_{s+2},\dots,K_l$ denote all the other $2$-handles of $X$. Since $N$ is simply connected, we may assume that each $K_j$ $(s+1\leq j\leq l)$ does not algebraically go over any $1$-handles of $N$, if necessary, by sliding the $K_j$ over $K_0,K_1,\dots,K_s$. Moreover, if necessary by sliding and isotoping handles $K_{s+1}, K_{s+2},\dots,K_l$, we may assume that $X-\{K_0,K_1,\dots,K_s\}$ is a good Legendrian handlebody and that $X$ is a Legendrian handlebody (cf.\ Section 4.2 of \cite{OS1}). Note that we keep handles of $N$ unchanged yet. Hence this Legendrian handlebody $X$ still contains the Stein handlebody $N$ as a subhandlebody.

Now we slide and isotope handles $K_0,K_1,\dots,K_s$ of $N$ so that $X$ becomes a good Legendrian handlebody, where we do not slide any $K_i$ $(0\leq i\leq s)$ over the handles $K_{s+1}, K_{s+2},\dots,K_l$. This finishes Step (I). For a while, we call this sliding operation of $K_0,K_1,\dots,K_s$ the ``normalization of $N$''. We use the same symbol $K_0,K_1,\dots,K_s$ for the resulting $2$-handles. We denote the resulting good Legendrian handlebody by $X^{(g)}$. 
Note that $X^{(g)}$ is diffeomorphic to $X$ and that we can recover the Stein handlebody decomposition of $N$ by reversing the operation ``normalization of $N$''.  

We next apply Steps (II) and (III) to the above $X^{(g)}$. Let $X_0,X_1,\dots,X_n$ be the resulting Legendrian handlebodies. It follows from Theorem~\ref{thm:AY5} that these $X_i$'s satisfy the claims (1)--(5) of Theorem~\ref{sec:thm:stein and non-stein}. Unfortunately, any Stein 4-manifold $X_i$ $(1\leq i\leq n)$ possibly do not contain the nucleus $N$. We thus proceed with this construction. 

We here apply cork twists to $X_n$ as follows (see also Figures~\ref{fig:nuclei_fig15}, \ref{fig:nuclei_fig16} and \ref{fig:nuclei_fig17}). Replace all the $W^+$-modifications applied to each $K_j$ $(0\leq j\leq s)$ of the Legendrian handlebody $X_n$ in Step (III) with the corresponding $W^-$-modifications. We denote by $\widetilde{X}_n$ the resulting Legendrian handlebody. Since $\widetilde{X}_n$ is obtained from $X_n$ by cork twists, $\widetilde{X}_n$ is homeomorphic to $X_n$. 

Let $\widetilde{X}$ be the Legendrian handlebody obtained from $X_n$ by removing all the $W_1$'s which links with $2$-handles of $N$. Namely, we define $\widetilde{X}$ to be the Legendrian handlebody obtained from $X_n$ by the following operations (i) and (ii) (see also Figures~\ref{fig:nuclei_fig15}, \ref{fig:nuclei_fig16} and \ref{fig:nuclei_fig17}). 
\begin{enumerate}
 \item [(i)] Remove all the $1$- and $2$-handles of $W_1$'s given by $W(p_i)$-modifications $(1\leq i\leq n)$ in Steps (II) and (III). 
 \item [(ii)] Remove all the $1$- and $2$-handles of $W_1$'s given by $W^+(q_j)$-modifications $(1\leq j\leq s)$ in Step (III). 
\end{enumerate}
As a smooth handlebody, $\widetilde{X}_n$ can be obtained from $\widetilde{X}$ by $W$-modifications. We can thus embed $\widetilde{X}_n$ into $\widetilde{X}$ by Lemma~\ref{lem:W:embedding}. Reversing the operation ``normalization of $N$'', we can change $\widetilde{X}$ into a Stein handlebody keeping its diffeomorphism type. Note that $N\subset \widetilde{X}_n \subset \widetilde{X}$.

Since $\widetilde{X}$ is a Stein handlebody, we can embed $\widetilde{X}$ (and hence $\widetilde{X}_n$) into a closed symplectic $4$-manifold with $b_2^+>1$. Theorem~\ref{sec:intro:thm:original} thus shows that $\widetilde{X}_n$ admits infinitely many distinct smooth structures, because $\widetilde{X}_n$ contains the nucleus $N$. Since $\widetilde{X}_n$ is homeomorphic to ${X}_n$, the rest of the proof follows from Theorem~\ref{sec:nonstein:thm:nonstein}. 
\end{proof}

Theorem~\ref{sec:thm:stein and non-stein} immediately gives the following. 

\begin{corollary}\label{cor:nonstein:2-handlebody:boundary sum}Let $Y$ be a $2$-handlebody, and let $N$ be a nucleus which admits a Stein structure. Put $X=Y\natural N$. Then, for each integer $n\geq 1$, there exist infinitely many compact connected oriented smooth $4$-manifolds $X_i$ $(i=0,1,2,\dots,)$ which satisfies $(1)$--$(5)$ of Theorem~\ref{sec:thm:stein and non-stein}.
\end{corollary}

Since every finitely presented group (resp.\ integral symmetric bilinear form $\mathbb{Z}^k\times \mathbb{Z}^k\to \mathbb{Z}$ $(k\geq 0)$) is realized  as the fundamental group (resp.\ intersection form) of a $2$-handlebody, this corollary shows the following. 

\begin{corollary}\label{cor:nonstein:group}$(1)$ Let $G$ be any finitely presented group. Then there exist arbitrary many compact Stein $4$-manifolds and infinitely many non-Stein smooth $4$-manifolds such that they are all homeomorphic but mutually non-diffeomorphic and that each of their fundamental groups is isomorphic to $G$.\\
$(2)$ Let $Q: \mathbb{Z}^k\times \mathbb{Z}^k\to \mathbb{Z}$ $(k\geq 0)$ be any integral symmetric bilinear form, and let $R: \mathbb{Z}^2\times \mathbb{Z}^2 \to \mathbb{Z}$ be any integral, indefinite, unimodular symmetric bilinear form. Then there exist arbitrary many simply connected compact Stein $4$-manifolds and infinitely many non-Stein smooth $4$-manifolds such that they are all homeomorphic but mutually non-diffeomorphic and that each of their intersection forms is isomorphic to $Q\oplus R$.
\end{corollary}




\end{document}